\documentclass[11pt]{article}
\usepackage{fullpage}
\usepackage{amsmath,amssymb,amsthm}
\usepackage{xcolor}
\usepackage{dsfont}
\usepackage{graphicx}
\usepackage{float}
\usepackage{subfig} 
\usepackage{stmaryrd}
\usepackage[style=alphabetic,giveninits=true]{biblatex}

\renewbibmacro{in:}{}
\definecolor{darkred}{RGB}{200,0,0}
\definecolor{darkgreen}{RGB}{0,150,0}
\definecolor{darkblue}{RGB}{0,0,200}
\usepackage[colorlinks=true,    
            linkcolor=darkred,  
            citecolor=darkgreen,    
            urlcolor=darkblue,      
            linktocpage=true,
            bookmarks=true]{hyperref}

\usepackage{thmtools}
\declaretheorem[numberwithin=section,name=Lemma,
refname={Lemma,Lemmas},
Refname={Lemma,Lemmas}]{lemma}
\numberwithin{equation}{section}

\date{}
\title{Numerical methods for the computation of densities of states of
periodic operators}
\author{Ewen Lallinec$^1$, Antoine Levitt$^1$}
\begin{document}
\maketitle

\begin{abstract}
	We present a comparative study of numerical methods for computing
	electronic densities of states (DOS) in periodic systems. We provide
	a detailed analysis of the domain of validity of the Brillouin
	complex deformation (BCD), a recently-proposed method promising
	exponential convergence without need for smearing. We compare on a
	range of systems the BCD with several methods, including the
	standard smearing and linear tetrahedron methods, as well as an
	adaptive integration method. Our results establish clear performance
	regimes for each method, offering practical guidance for DOS
	computations across a range of systems and accuracy
	requirements.
\end{abstract}
\footnotetext[1]{Université Paris-Saclay, CNRS, Laboratoire de mathématiques d’Orsay, 91405, Orsay, France}

\section{Introduction}\label{section:Introduction}

A $d$-dimensional perfect
crystal is modeled using a Hamiltonian $H$, periodic on a lattice
$\mathcal{R}\subset \mathbb{R}^d$. We consider energy-resolved properties of the system of
the form $I(E) = {\underline{{\rm Tr}}} (F \delta(H-E))$, where ${\underline{{\rm Tr}}}$ is
the trace per unit volume and $F$ is a periodic observable. Because of
periodicity, we can express this as an integral over the Brillouin
zone $\mathcal{B}$, a unit cell of the reciprocal lattice
\begin{align}
	I(E)= \frac 1 {|\mathcal{B}|}\int_{\mathcal{B}}\sum_{n}f_{n \mathbf{k}}\mathds{\delta}(\varepsilon_{n\mathbf{k}}-E)\mathrm{d}\mathbf{k}\label{eq:bzintegral},
\end{align}
with $\varepsilon_{n\mathbf{k}}$ the energy bands (eigenvalues of
the Bloch-transformed Hamiltonian $H_\mathbf{k}$), $E$ an energy parameter and $f_{n\mathbf{k}}$ a
function of $\mathbf{k}$ (the matrix element
$\langle  u_{nk}|F_{k}|u_{nk} \rangle$, where $u_{nk}$ is the periodic part of the
Bloch wave $\psi_{nk}$). The density of states (DOS) is a specific case of the
integral~\eqref{eq:bzintegral} where $F=1$:
\begin{equation}
	D(E) = \frac 1 {|\mathcal{B}|}\int_\mathcal{B}\sum_{n} \mathds{\delta}(\varepsilon_{n\mathbf{k}}-E)\mathrm{d}\mathbf{k}.
\end{equation}
Physically, $D(E) dE$ quantifies the number of electronic states per
unit volume available in an infinitesimal energy range $[E,E+dE]$.
This makes it a crucial quantity for understanding and predicting bulk
properties of materials. Although this study focuses on the
computation of the DOS, the algorithms we study can be extended to
more general Brillouin zone integration problems such as Green
functions or conductivities~\cite{duchemin2023efficient, vanmunoz2025highorder}.

The integrand contains singularities which are concentrated on the
sheets of the Fermi surface
$S_{n}(E)=\left\{ \mathbf{k} \in \mathcal{B},\ \varepsilon_{n\mathbf{k}}=E \right\}$~\cite{cances2020numerical}.
The geometry of this surface can pose significant computational
issues. To compute these integrals, a common approach involves the
computation of a smoothed DOS which depends on a smearing parameter
$\eta$. This smeared DOS is then non-singular and is defined as
\begin{equation}
	D_{\eta}(E) = \frac 1 {|\mathcal{B}|}\int_{\mathcal{B}} \sum_{n}\mathds{\delta}_\eta(\varepsilon_{n\mathbf{k}}-E)\mathrm{d}\mathbf{k},{\rightarrow} D(E),
\end{equation}
with $\delta_\eta$ a smooth kernel approximating the Dirac distribution in the
limit $\eta \rightarrow 0^+$, such that $D_{\eta} \to D$ as $\eta \to 0$ in the sense of distributions. The kernel can commonly be chosen to be
Gaussian, Lorentzian, or equal to minus the derivative of the
Fermi-Dirac distribution. Higher-order kernels can also be chosen
\cite{methfessel1989highprecision,cances2020numerical,colbrook2021computing}
to achieve faster asymptotic convergence as $\eta\to0$, at the price
of some artifacts (negative or oscillatory densities of states), but
we will focus in this paper on the simple Lorentzian kernel.

Once smeared, the integrand can be discretized using the standard
periodic trapezoidal rule (PTR \cite{trefethen2014exponentially},
also referred to as Monkhorst-Pack grid \cite{monkhorst1976special}), which
converges exponentially for analytic and periodic integrands.
Recently, adaptive methods have been proposed, such as the iterative
adaptive integration method (IAI)~\cite{kaye2023automatic}, promising
faster convergence for small $\eta$, or when the integrand is
regularized by self-energies with small anti-hermitian parts (such as
in the DMFT method).

Another class of methods targets the non-smeared integral
\eqref{eq:bzintegral} directly. The standard method is the Linear
Tetrahedron method
(LT)~\cite{lehmann1972numerical,jepson1971electronic,blochl1994improved},
which divides the Brillouin zone into regular tetrahedra and computes
the DOS by linearly interpolating the bands and matrix elements within
each tetrahedron and integrating analytically the DOS of this linear
interpolant. A more recent method is the Brillouin Complex Deformation
(BCD), which is based on a complex deformation of the Brillouin zone
to avoid singularities. This method has been recently introduced by one of us
\cite{duchemin2023efficient} for the numerical computation of
resonances, and is similar to other ideas used in scattering theory
\cite{gerard1990resonance,gerard1998mourre} or
in the far-field asymptotics of Green functions \cite{assier2023contribution,shanin2024double}.

In this article, we compare these four different approaches on a
common footing, with the aim of assessing their accuracy, efficiency,
and range of applicability for the computation of densities of states.
The regime that we study in this paper is that of accurate results on
zero-smearing (or low-smearing) systems. Particular attention is paid
to the Brillouin complex deformation method, for which we provide a
more detailed discussion, including a careful analysis of its domain
of validity in the presence of band crossings.

We benchmarked the methods on a set of representative systems of
increasing complexity, ranging from simple one-dimensional models to
realistic two- and three-dimensional Hamiltonians. For each system, we
consider several representative energies, chosen to span different
levels of numerical difficulty (from smooth regions to energies close
to van Hove singularities), in order to cover the various regimes
encountered in practical DOS computations. This allows us to highlight
the regimes in which smearing-based approaches are appropriate, as
well as those in which smearing-free methods, and in particular the BCD,
provide a significant advantage.

\section{Theoretical settings}\label{section:framework}
\subsection{Tight-binding Hamiltonian and band structure}
For concreteness, we consider a tight-binding model. The methods
discussed are in principle applicable to continuous models with infinitely many bands, but since the integrals
we consider are limited to a finite number of bands, it is more
advantageous computationally to consider a Wannier interpolation \cite{yates2007spectral}, which reduces to a tight-binding model.

Let $\mathcal{R} \cong \mathbb{Z}^d$ be a lattice on $\mathbb{R}^d$ with
$d$ the dimension of our physical space. The Hamiltonian
$H$ acts on wave functions $\psi$ living in the Hilbert space
$\mathcal{H}\cong \ell^2(\mathbb{Z}^d,\mathbb{C}^M)$, where $M$
represents the number of degrees of freedom per unit cell. Consequently,
the Hamiltonian $H(\mathbf{R},\mathbf{R'})$ between two lattice sites
$\mathbf{R}$ and $\mathbf{R'} \in \mathbb Z^{d}$  is represented as a
matrix in $\mathbb{C}^{M\times M}$. It satisfies the periodicity condition for $\mathbf{T} \in \mathbb{Z}^d$
\begin{equation}
	H(\mathbf{R}+\mathbf{T},\mathbf{R'}+\mathbf{T}) = H(\mathbf{R},\mathbf{R'}).
\end{equation}
The Brillouin zone is defined as $\mathcal{B}\cong\left[-\pi,\pi\right)^d$. Thanks to $H$
periodicity, its generalized eigenstates can be expressed as Bloch waves $\psi_{n\mathbf{k}}$ (with $\mathbf{k} \in \mathcal{B}$) satisfying
\begin{equation}
	\psi_{n\mathbf{k}}(\mathbf{R}) =u_{n\mathbf{k}}e^{i\mathbf{k}\cdot\mathbf{R}},
\end{equation}
where $u_{n\mathbf{k}} \in \mathbb{C}^M$ are $\mathcal{R}$-periodic, orthonormal eigenvectors of the reciprocal space Hamiltonian $H_\mathbf{k}$:
\begin{equation}
	H_{\mathbf{k}}u_{n\mathbf{k}}=\varepsilon_{n\mathbf{k}}u_{n\mathbf{k}}.
	\label{eq:schrodinger}
\end{equation}
The reciprocal space Hamiltonian $H_{\mathbf{k}}$ is obtained as the Bloch (Fourier) transform of $H(\mathbf{0},\cdot)$ at point $\mathbf{k}$ of $\mathcal{B}$
\begin{equation}
	H_{\mathbf{k}} = \sum_{\mathbf{T} \in \mathcal{R}}e^{i\mathbf{k}\cdot\mathbf{T}}H(\mathbf{0},\mathbf{T}).
\end{equation}
In practice, $H(\mathbf{R},\mathbf{R'})$ has a finite range, and
therefore $H_{\mathbf{k}}$ is an entire function of
$\mathbf k \in \mathbb C^{d}$ (a trigonometric polynomial).

\subsection{Density of states}
\label{sec:dos}
The density of states (DOS) at energy $E\in \mathbb{R}$ is formally defined (in the sense of distributions) as
\begin{equation}
	D(E)  =\frac{1}{\left|\mathcal{B}\right|}\int_{\mathcal{B}}\sum_{n=1}^M\delta(E-\varepsilon_{n\mathbf{k}})\mathrm{d}\mathbf{k}.
	\label{eq:dosdirac}
\end{equation}
A change of variable also yields~\cite{cances2020numerical}

\begin{equation}
	\ D(E) = \frac{1}{\left|\mathcal{B}\right|} \sum_{n=1}^M \int_{S_n(E)} \frac{1}{|\nabla\varepsilon_{n\mathbf{k}}|} \mathrm{d}\mathcal{H}^{d-1}(\mathbf{k}).
	\label{eq:dosgrad}
\end{equation}
with $\mathrm{d}\mathcal{H}^{d-1}$ the Hausdorff measure and $S_n(E)$ and $S(E)$ represent the isosurfaces defined by

\begin{align}
	S_n(E) & = \left\{k \in \mathcal{B}, \ \varepsilon_{n\mathbf{k}}=E\right\}, \\
	S(E)   & = \bigcup_{n=1}^M\mathcal{S}_n(E).
	\label{eq:isosurfaces}
\end{align}
The DOS can be shown to be smooth at $E$ provided that~\cite{cances2020numerical}
\begin{itemize}
	\item The bands do not cross at energy $E$ \begin{equation}\forall \ n, m \in \llbracket 1,M \rrbracket, n\neq m,  S_{n}(E) \cap S_{m}(E) = \emptyset.
		      \label{eq:cond1}\end{equation}
	\item The band gradients (Fermi velocities) are non-zero on the Fermi surface \begin{equation}\forall n \in \llbracket 1,M \rrbracket, \mathbf{k}\in S_{n}(E), \ \nabla \varepsilon_{n\mathbf{k}} \neq 0. \label{eq:cond2}\end{equation}
\end{itemize}

Under these conditions, the Fermi surface $S(E)$ is locally defined by
a single condition $\varepsilon_{n\mathbf k} = E$ with $\nabla \varepsilon_{n\mathbf k} \neq 0$ and is
therefore a smooth surface. The energies at which the DOS is not
smooth (caused by the violation of one or both of these conditions)
are called van Hove singularities; however, there are cases where the
first condition is violated at non-van-Hove energies (see~\autoref{section:crossings} for a discussion).

\section{Smearing methods}\label{section:smearing}
The smeared DOS with kernel $\delta_{\eta}(x) = \eta^{-1} \delta_{1}( x/ \eta)$, an approximation of the
delta function, is given by
\begin{align}
	D_{\eta}(E) & = \frac 1 {|\mathcal{B}|}\int_{\mathcal{B}} \sum_{n}{\delta}_\eta(\varepsilon_\mathbf{k}-E)\mathrm{d}\mathbf{k}, \\
	            & = \int_{-\infty}^{+\infty} \delta_{\eta}(E-E') D(E')dE',                                                         \\
	            & = \int_{-\infty}^{+\infty} \delta_{1}(E')D(E+\eta E') dE'.
	\label{eq:smearing_error}
\end{align}
Expanding in powers of $\eta$, we see formally that the error
$D_{\eta}(E) - D(E)$ is of
order $\eta^{2}$ when $D$ is smooth at $E$. Indeed this is rigorous
for kernels with sufficient decay at infinity
\cite{cances2020numerical}, but because the Lorentzian kernel has
heavy tails (no finite second moment) this argument is invalid.
However, identifying the smeared density of states with the
Green function\eqref{eq:greendos} and using the analyticity
properties of the Green function (which can be shown
e.g., using the BCD method) shows the expansion is indeed valid.

Other popular kernels include the gaussian or Fermi-Dirac
kernels, but numerical tests revealed no significant difference
among these three kernels so that we focus on the
Lorentzian kernel in this paper. Higher-order kernels
\cite{methfessel1989highprecision,cances2020numerical,
	colbrook2021computing} are asymptotically more efficient but
are not investigated in this work.

Once regularized, the approximate DOS $D_{\eta}$ can be evaluated by
using a quadrature method, yielding a combined approximation
$D_{N,\eta}$ such that
\begin{equation}
	D(E) = \lim_{\eta\rightarrow 0^+}\lim_{N\rightarrow +\infty} D_{N,\eta}(E).
\end{equation}
where the interchange of limits is invalid because of the singularity
of the integrand as $\eta \rightarrow 0^+$.

Even after regularization, at a fixed but small $\eta$, computing $D_{\eta}(E)$ is not a simple task
because the integrand will be sharply peaked around the Fermi
surface $S(E)$. More precisely, by a Taylor expansion we see
that the integrand is a bump function located near the Fermi
surface with a width orthogonal to the Fermi surface equal to
$\eta/|\nabla \varepsilon_{n\mathbf{k}}|$.

\subsection{Periodic Trapezoidal Rule}

To compute the smeared DOS $D_\eta$, we discretize the Brillouin zone
$\mathcal{B}$ using the Monkhorst-Pack grid~\cite{monkhorst1976special}, a $N \times N \times N$ (in
3D) uniform grid of $\mathcal B$ defined by
\begin{equation}
	\mathcal{B}_N = \left\{\frac{2\pi}{N}\mathbf{n}-\pi, 0\leq n_i \leq N-1,i \in \left\{1,\dots,d\right\}\right\}.
\end{equation}
yielding the approximation
\begin{equation}
	\frac{1}{|\mathcal{B}|}\int_\mathcal{B}f(\mathbf{k})\mathrm{d}\mathbf{k} \approx \frac{1}{N^d}\sum_{\mathbf{k} \in\mathcal{B}_N} f(\mathbf{k}).
	\label{eq:Monkhorstregular}
\end{equation}

The previous discussion can be made more precise: the rate of
convergence of the PTR method depends on the decay of the
Fourier coefficients of the integrand
\cite{trefethen2014exponentially}. For an analytic integrand, the
Fourier coefficients decay exponentially at a rate equal to the
distance to the singularity of $f$ in $\mathbb C^{d}$ that is closest to the real
space $\mathbb R^{d}$. Focusing on a single band in 1D
for simplicity, the integrand is $1/(E+i\eta-\varepsilon_{\mathbf{k}})$, which is
analytic as long as its denominator is non-zero. Taylor
expanding near a point $\mathbf{k}_{0}$ of the Fermi surface, we find that the
denominator vanishes for $\mathbf k \in \mathbb C$ of imaginary part
close to $\eta/\varepsilon_{n\mathbf k_{0}}'$ for $\eta$ small. Therefore, the
error behaves asymptotically as $e^{-\alpha \eta L}$ for some $\alpha > 0$. In
particular, the convergence is very slow for small $\eta$.

\subsection{Iterated Adaptive Integration}\label{IAI}
The Iterated Adaptive Integration (IAI)~\cite{kaye2023automatic} was developed to compute the smeared DOS for small smearing parameters $\eta$. Since smearing introduces localized peaks in the DOS, an adaptive discretization of the Brillouin zone is bound to be more efficient than uniform sampling.

The adaptive Gaussian quadrature used in the IAI algorithm is recursive over the
dimension, and adaptively decomposes the Brillouin zone (considered as
a cube) into smaller rectangles. This adaptive meshing efficiently
captures the localized peaks in the smeared DOS. Consequently, the IAI method
achieves a better asymptotic dependence on $\eta$ than the PTR, with a
convergence rate that degrades only logarithmically with $\eta$
\cite{kaye2023automatic} (as opposed to linearly for the PTR). It is therefore
able to achieve high precision even at small $\eta$. Its main
disadvantage is its high computational overhead.

\section{Smearing-free methods}\label{section:methods}
We now discuss methods for computing the density of states directly,
i.e. without smearing.
\subsection{Linear Tetrahedron}\label{LT}

The Linear Tetrahedron method
(LT)~\cite{lehmann1972numerical,jepson1971electronic,blochl1994improved}
is a widely used approach for Brillouin zone integration. It relies on
linearly interpolating the bands within a tetrahedral decomposition of
the Brillouin zone. To achieve this, the Brillouin Zone is first
divided into $N^d$ cubes, each of which is further subdivided into six
tetrahedra. Within each tetrahedron, the bands are linearly
interpolated to approximate the DOS. This approach yields an explicit
analytical expression (directly computed from~\autoref{eq:dosgrad}) for
the DOS contribution in each tetrahedron as a function of energy $E$.

The LT method achieves a relatively slow convergence rate of
$\mathcal{O}\left( N^{-1}\right)$. This can be most easily seen in 1D,
where the LT method approximates the values of $1/|\varepsilon'_{k}|$
on the Fermi surface by finite differences, therefore being no more
than $1/N$ accurate. The LT method is appealing by its conceptual
simplicity (it does not depend on any adjustable parameters) and is
often very useful in low-accuracy requirements, but its slow
convergence makes it hard to obtain precise results. Its convergence
can also be somewhat erratic, because small band
gradients can yield divergent results for the LT method near van Hove
singularities (see \cite{kaye2023automatic} Appendix E for a
discussion). Finally, its implementation presents difficulties in the
case where a self-energy is used (as is the case in DMFT for
instance) \cite{kaye2023automatic}.

\subsection{Brillouin Complex Deformation}\label{BCD}
\paragraph{The Green function}
We specialize the discussion of the smeared DOS to the Lorentzian kernel
\begin{align*}
	K_{\eta}(E) = \frac 1 \pi \frac \eta {\eta^{2}+E^{2}} = - \frac 1 \pi {\rm Im} \frac 1 {E+i\eta}.
\end{align*}
This last expression reveals an interesting connection between the
smeared DOS
\begin{align}
	D_{\eta}(E) & = \frac 1 {|\mathcal{B}|}\int_{\mathcal{B}} \sum_{n=1}^M{K}_\eta(\varepsilon_\mathbf{k}-E)\mathrm{d}\mathbf{k}, \\
	            & =  \int_{-\infty}^{+\infty} K_{\eta}(E-E') D(E')dE',
	\label{eq:smearing_error2}
\end{align}
and the Green function, defined by
\begin{align}
	G(\mathbf{R},\mathbf{R}';E+i\eta) & =  \frac{1}{\left|\mathcal{B}\right|} \int_{\mathcal{B}}e^{i\mathbf{k}\cdot(\mathbf{R}-\mathbf{R}')}\frac{1}{E+i\eta-H_{\mathbf{k}}}\mathrm{d}\mathbf{k},                                                    \\
	                                  & = \sum_{n=1}^M\frac{1}{\left|\mathcal{B}\right|}\int_{\mathcal{B}}e^{i\mathbf{k}\cdot(\mathbf{R}-\mathbf{R}')}\frac{u_{n\mathbf{k}}u_{n\mathbf{k}}*}{E+i\eta-\varepsilon_{n\mathbf{k}}}\mathrm{d}\mathbf{k}.
	\label{eq:green}
\end{align}
Using the identity
$\delta(x) = -1/ \pi \, {\rm Im}\lim_{\eta \to 0} 1 /(x+i\eta)$ (in the sense of distributions), we get
\begin{align}
	D(E) & = -\frac{1}{\pi}\lim_{\eta\rightarrow 0^+}\mathrm{Im}\left(\underline{\mathrm{Tr}}(G(\mathbf{0},\mathbf{0};E+i\eta))\right),                                                                        \\
	     & = -\frac{1}{\pi}\lim_{\eta\rightarrow 0^+}\mathrm{Im}\left(\frac{1}{\left|\mathcal{B}\right|}\int_{\mathcal{B}}\mathrm{Tr}\left(\frac{1}{E+i\eta-H_{\mathbf{k}}}\right)\mathrm{d}\mathbf{k}\right).
	\label{eq:greendos}
\end{align}
\paragraph{Complex deformation}
The Brillouin Complex Deformation (BCD)~\cite{duchemin2023efficient}
is a method originally developed to compute resonances, which we now
apply to the calculation of the density of states. The basic idea of
this method is to deform the Brillouin zone into the ($\mathbb{C}^{d}$) complex
domain, in order to avoid singularities of the integrand in \autoref{eq:greendos}. The
resulting integral has a smooth integrand, and is then discretized
using the standard PTR method. More precisely, we use the following
easy lemma, which can be viewed as a complex change of variable, or a
generalization of Cauchy's formula to multidimensional integrals \cite{duchemin2023efficient}:
\begin{lemma}
	Let $I$ be a $\mathcal{R}^*$-periodic function, analytic in an open set $U=\mathbb{R}^d+i[-\eta,\eta]^d$. Then, for all $\mathcal{R}^*$-periodic and continuously differentiable functions $\mathbf{h}(\mathbf{k}):\mathbb{R}^d\rightarrow[-\eta,\eta]^d$, we have

	\[\int_{\mathcal{B}}I(\mathbf{k})\mathrm{d}\mathbf{k}=\int_{\mathcal{B}}I(\mathbf{k}+i\mathbf{h}(\mathbf{k}))\det(\mathbf{1}+i\mathbf{h}'(\mathbf{k}))\mathrm{d}\mathbf{k}.  \]\label{lemma:BCD}
\end{lemma}
Our goal is to choose a deformation $\mathbf{h}(\mathbf k)$ such
that the spectrum of $H_{\mathbf{k}+i\mathbf{h}(\mathbf{k})}$
near $E$ gets pushed down in the complex
plane, so that $1/(E+i\eta-H_{\mathbf{k}+i\mathbf{h}(\mathbf{k})})$
stays analytic as a function of $\mathbf k$. If such a deformation can be found, then we have

\begin{equation}
	D(E)=-\frac{1}{\pi}\mathrm{Im}\left( \frac{1}{\left|\mathcal{B}\right|}\int_{\mathcal{B}} \mathrm{Tr}\left(\frac{1}{E-H_{\mathbf{k}+i\mathbf{h}(\mathbf{k})}}\right)\det\left( \mathbf{1}+i\mathbf{h'}(\mathbf{k}) \right)\mathrm{d}\mathbf{k} \right).\label{eq:bcddos}
\end{equation}
and the integrand is now non-singular.

\paragraph{Choosing the deformation}
To choose a suitable deformation $\mathbf h$, consider $E\in \mathbb{R}$ such
that there exists $n$ and $\mathbf{k}_0$ with $\varepsilon_{n \mathbf{k}_0}=E$,
with $\varepsilon_{n\mathbf k}$ smooth near $\mathbf k_{0}$. Near $\mathbf{k}_0$, the band can be approximated at first-order by
\begin{equation}
	\varepsilon_{n\mathbf{k}}\approx E+\nabla\varepsilon_{n\mathbf{k}_0}\cdot(\mathbf{k}-\mathbf{k}_0).
	\label{eq:linearization}
\end{equation}
After applying the BCD, the imaginary part of the band becomes
\begin{equation}
	\mathrm{Im}(\varepsilon_{n\mathbf{k}+i\mathbf{h}(\mathbf{k})})\approx \nabla\varepsilon_{n\mathbf{k}_0}\cdot\mathbf{h}(\mathbf{k}).
\end{equation}
To ensure this is negative near $S(E)$, we therefore require
\begin{equation}
	\mathbf{h}(\mathbf{k})\cdot\nabla\varepsilon_{n\mathbf{k}}<0.
\end{equation}

To ensure this condition for a single band, we could simply take
$\mathbf{h}(\mathbf{k}) = -\alpha \nabla\varepsilon_{n\mathbf{k}}$ for some $\alpha > 0$. However, for
multiple bands, we need to ensure the smoothness of
$\mathbf{h}$. We achieve this by selecting the bands that are
close to the Fermi surface \cite{duchemin2023efficient}:
\begin{align}
	\mathbf{h}(\mathbf{k}) & =-\alpha\sum_{n=1}^M\nabla_\mathbf{k}\varepsilon_{n\mathbf{k}}\chi\left(\frac{\varepsilon_{n\mathbf{k}}-E}{\Delta E}\right),   \\
	                       & =-\alpha \left(\mathrm{Tr}\left(\nabla_{\mathbf k} H_\mathbf{k}\chi\left(\frac{H_\mathbf{k}-E}{\Delta E}\right)\right)\right),
	\label{eq:deformation}
\end{align}
where
$\chi$ is a smooth cutoff function, typically a Gaussian $\chi(x) = e^{-x^{2}}$. In
particular, the second formulation in~\autoref{eq:deformation} shows the
smoothness of $\mathbf{h}$. The computation of the derivatives of $\mathbf h(\mathbf k)$ can be
performed using standard perturbation theory. We detail these issues
in~\autoref{divdif}.

When $E$ satisfies the two conditions in~\autoref{sec:dos} (no band
crossing at the Fermi level, and non-zero band gradients at the Fermi
level), and when $\Delta E$ and $\alpha$ are small, the linearization
\eqref{eq:linearization} is valid and the BCD scheme is guaranteed to
regularize the integral. In practice however, the non-crossing condition might
be too stringent, as we will see in~\autoref{section:crossings}.

\begin{figure}[ht!]
	\centering
	\subfloat{\includegraphics[scale=0.24]{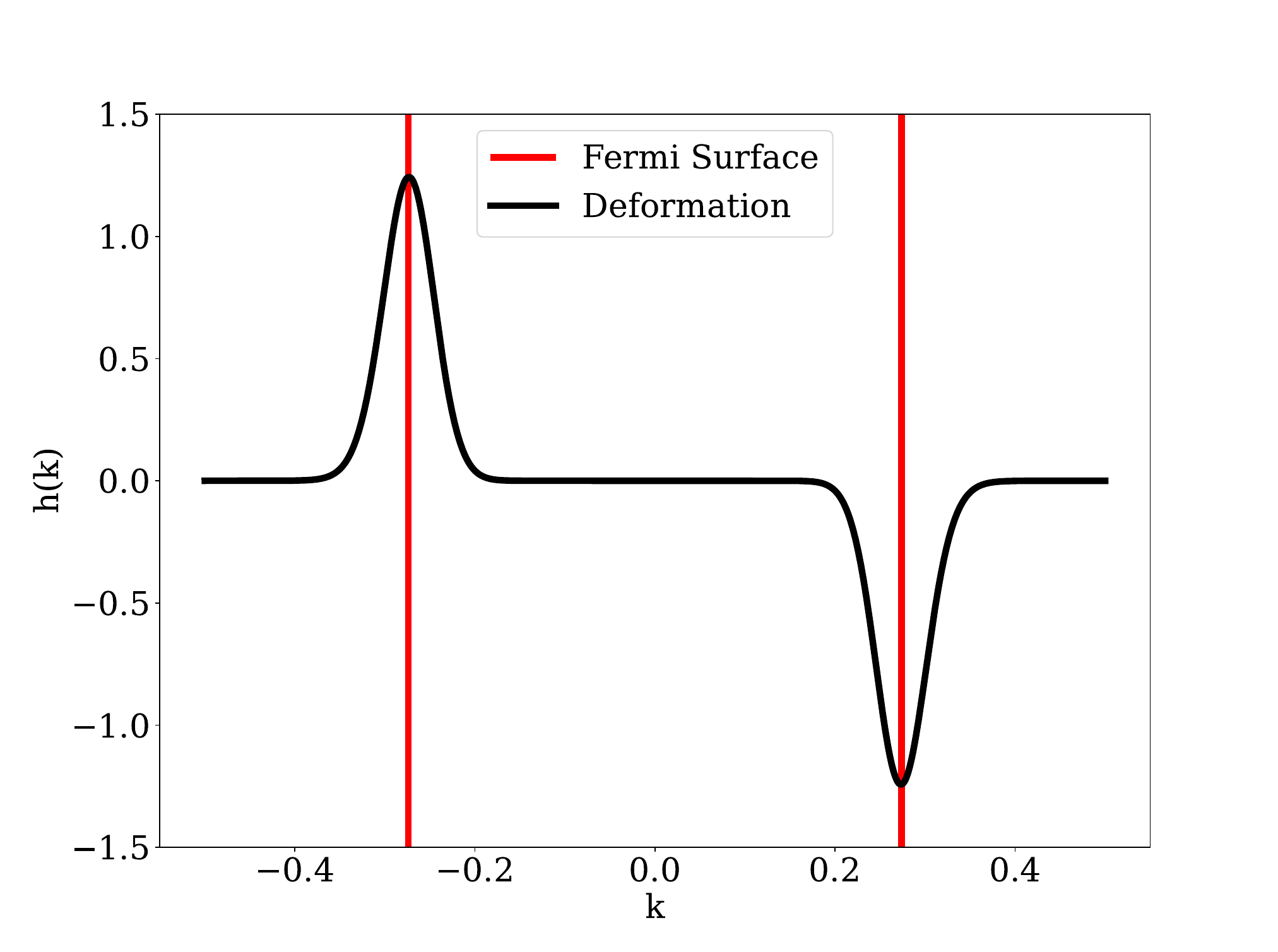}}\hfill
	\subfloat{\includegraphics[scale=0.24]{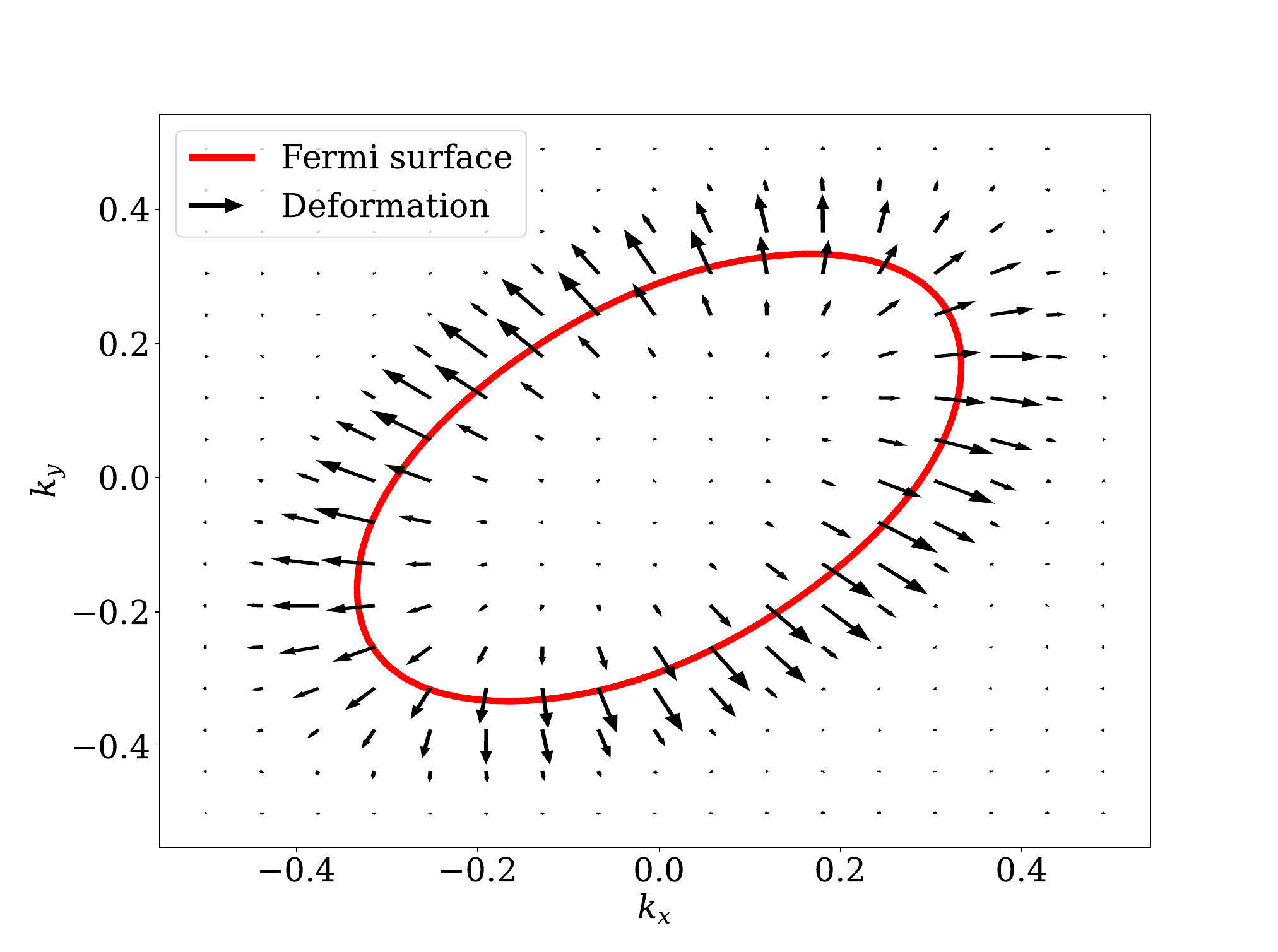}}
	\caption{Examples of deformation for the monatomic chain in 1D
		(on the left) and graphene in 2D (on the right). We represent as vectors the complex deformation in 2D.}
	\label{fig:exampledef}
\end{figure}

To illustrate the deformation process, we provide examples of $\mathbf{h}(\mathbf{k})$ in 1D and 2D for simple isosurfaces $S(E)$ in~\autoref{fig:exampledef}.
\paragraph{Choosing the parameters}

The deformation function $\mathbf{h}$ depends on the parameters
$\alpha$ and $\Delta E$. In principle, under the conditions of \autoref{sec:dos}, the scheme is guaranteed to successfully regularize the
integral when $\Delta E$ and $\alpha$ are sufficiently small, but taking these
parameters too small results in sharply peaked integrands.

The parameter $\Delta E$ controls the energy width near the Fermi level around which band
gradients are taken into account. If selected too small, this would
result in a sharply varying $\mathbf h(\mathbf k)$, which would be
numerically difficult to integrate. If selected too large, the
contribution of the ``correct'' band gradient
$\nabla \varepsilon_{n \mathbf k}$ can be overwhelmed by that of other
bands, resulting in $\mathbf h$ pointing in the wrong direction
relative to $\nabla \varepsilon_{n \mathbf k}$: this would push the
spectrum up instead of down, and result in the wrong density of states
being computed. This is likely to happen in regions where the density
of states itself has particularly sharp variations.

The parameter $\alpha$ controls the magnitude of the deformation. It should
be selected sufficiently large so that the deformation pushes the
spectrum significantly away from the real axis, but not so large that
the linear approximation becomes invalid. Assuming $\Delta E$ is
chosen appropriately so that only one band $\varepsilon_{\mathbf k}$
is relevant, we can expand in powers of $\alpha$
\begin{align*}
	\varepsilon_{\mathbf{k} + i \mathbf h(\mathbf k)} \approx E - i \alpha |\nabla\varepsilon_{\mathbf{k}}|^{2} - \frac 1 2 \alpha^{2} \nabla^{2} \varepsilon_{\mathbf k} [\nabla \varepsilon_{\mathbf k}, \nabla \varepsilon_{\mathbf k}].
\end{align*}
We see that stopping the expansion to first order is legitimate when
$\frac 1 2\alpha |\nabla^{2}\varepsilon_{\mathbf k}| \ll 1$, so that
$\alpha$ should be selected smaller than the relevant inverse curvatures in
the system. Note that this criterion is somewhat pessimistic (in particular, the second-order term
contributes only a real part and does not modify the imaginary part,
so the expansion could be pushed to third order) but simple.

In practice, we recommend inspecting the band structure (of the system
under consideration or of a similar system) and setting $\Delta E$ to
be of the order of magnitude of the variations of the density of
states, and $\alpha$ of the order of magnitude of the inverse
curvatures in the system. In practice, we did not particularly try to
optimize these settings and simply took $\Delta E \sim 0.3-0.1$ eV and
$\alpha \sim 0.1 \text{ eV}^{-1}$ for all the systems considered (in
relative coordinates where the Brillouin zone is a cube of length $1$).

\section{Validity of the BCD}\label{section:validity}
As we discussed, the BCD is valid for isolated bands, regularizing the
integrand in the absence of a van Hove singularity. The von
Neumann-Wigner theorem \cite{vonneumann1993uber} states that the hermitian
matrices with repeated eigenvalues, considered as a subset of the
manifold of hermitian matrices, can be seen as the intersection of
manifolds of codimension at least $3$. As a consequence, generically,
one does not expect any eigenvalue crossings in dimensions $1$ and
$2$, and in 3D one expects them at isolated points. Those points are
generically Weyl points (3D conical crossings) and give rise to van
Hove singularities. Therefore, \emph{generically}, one would expect the
BCD to be efficient at all non-van Hove energies, where we expect any method to struggle anyway because of the roughness of the DOS at these energies.

However, as we will see, non-Weyl
crossings do happen robustly in realistic Hamiltonians, in any
dimensions, due to the presence of (non-generic) symmetries. There, the BCD method is at best heuristic. We will show in specific examples that it in fact does fail in some cases.

\subsection{Two independent bands in 1D}
\label{section:crossings}

\begin{figure}[ht!]
	\centering
	\includegraphics[scale=0.29]{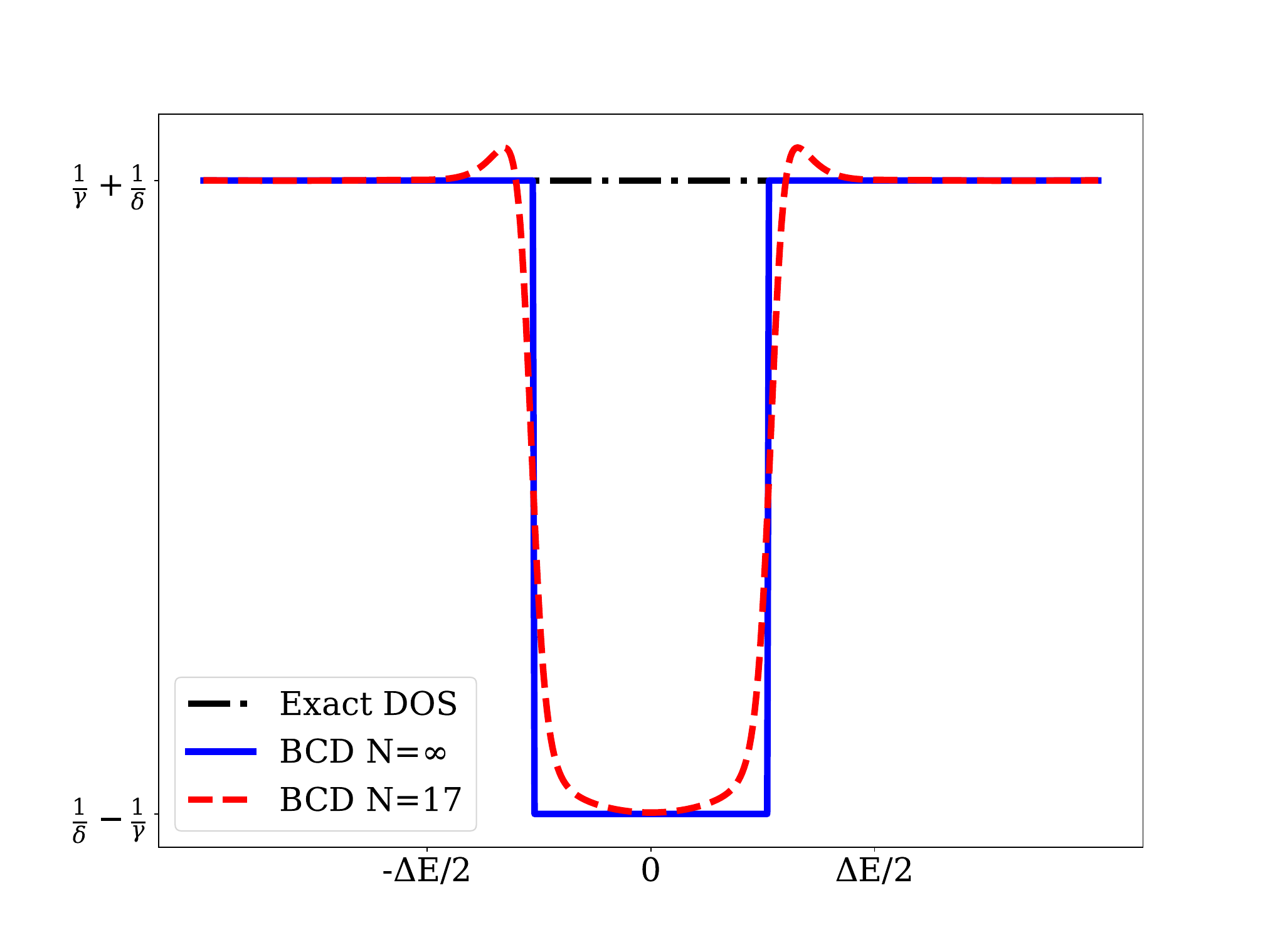}
	\caption{DOS computed with the BCD, displaying a systematic
		underestimation near $0$. With a finite sampling,
		oscillations occur near the transition.}
	\label{fig:toymodelbcd}
\end{figure}

Our choice of deformation in the BCD method was motivated by the single band case. To extend to the multiband case, we assumed that all bands at the same energy have a compatible gradient. A trivial case where this assumption is wrong is when the Hamiltonian is made up of two independent blocks. For instance, consider the toy model on $\mathcal{B}=\mathbb{R}$
\begin{equation}
	H(k)= \begin{pmatrix}
		\gamma k & 0         \\
		0        & -\delta k
	\end{pmatrix},
\end{equation}
where $ \gamma,\delta > 0$ and $\gamma\neq \delta$.
A straightforward calculation using~\autoref{eq:dosgrad} shows that the
DOS is (up to normalization)
\begin{equation}
	D(E) = \left( \frac{1}{\gamma} +\frac{1}{\delta}\right).
\end{equation}
Taking the deformation $h$ as defined for the BCD gives
\begin{equation}
	h(k) = - \alpha\gamma\exp\left( -\left(\frac{\gamma k-E}{\Delta E}\right)^2 \right) +\alpha\delta\exp\left( -\left(\frac{-\delta k-E}{\Delta E}\right)^2 \right).
\end{equation}
At $E=0$ and point $k=0$, this yields
\begin{equation}
	H(k+ih(k))= \begin{pmatrix}
		i\alpha \gamma(\delta-\gamma) & 0                             \\
		0                             & i\alpha \delta(\gamma-\delta)
	\end{pmatrix}.
\end{equation}
The diagonal components have opposite imaginary part. Thus one pole is
shifted into the upper complex plane and the other in the lower
complex plane, depending on the sign of $\gamma	-\delta$. The BCD will then count
\emph{negatively} the contribution of the pole shifted downwards, giving a
wrong result. This irregularity extends over an energy range $\sim \Delta E$,
due to the Gaussian factor scaling as $E/\Delta E$.
We illustrate this process in~\autoref{fig:toymodelbcd}.

As expected, a gap of width proportional to $\Delta E$ appears around $0$,
where the contribution of the second eigenvalue $1/\delta$ is negatively
added. Computing the BCD value of the DOS with a finite number $N$ of
points, as done in practice, introduces oscillations as
discontinuities (related to the Gibbs phenomenon). These oscillations
vanish as $N$ increases, but the gap remains.

\subsection{Crossings and symmetries}
\label{sec:crossings_symmetries}

The previous example was particularly simple but artificial. The fact
that the bands crossed is non-generic; almost any off-diagonal
perturbation will turn the crossing into an avoided crossing. We
investigate how the BCD behaves in the presence of general band
crossings, and how such band crossings can happen in more realistic systems.

Consider a
$2 \times 2$ Hamiltonian for simplicity. Such a Hamiltonian can be put
in the form
\begin{equation*}
	H(\mathbf{k}) = a(\mathbf{k}) + b(\mathbf{k}) \cdot \sigma,
\end{equation*}
where $b(\mathbf{k}) \in \mathbb R^{d}$ and $\sigma_{1}, \sigma_{2}, \sigma_{3}$ are the Pauli
matrices. At a point $\mathbf{k}_{0}$ of multiplicity $2$ at energy $E$,
assuming that $\chi$ is even (typically the gaussian kernel), the BCD prescription gives
\begin{equation*}
	\mathbf h(\mathbf k_{0}) = -\alpha {\rm Tr}\left(\nabla H_{\mathbf k_{0}} \chi\left(\frac{H_{\mathbf k_{0}}-E}{\Delta E}\right)\right) = -\alpha \nabla {\rm Tr}(H_{\mathbf k_{0}}) = -2 \alpha \nabla a(\mathbf k_{0}).
\end{equation*}
The deformed Hamiltonian is then
\begin{equation}
	H_{\mathbf{k} + i \mathbf h(\mathbf k)} \approx a(\mathbf{k}) -2 \alpha i \left( |\nabla a|^{2} + \sum_{i=1}^3 {(\nabla b_{i}\cdot \nabla a)} \sigma_{i}\right). \label{eq:bad_deformation}
\end{equation}
thus the eigenvalues are pushed off the real axis by an
amount $-2\alpha i \left(|\nabla a|^{2} \pm \sqrt{\sum_{i=1}^{3} |\nabla b_{i}\cdot \nabla a|^{2}}\right)$.
When
\begin{equation*}
	\sqrt{\sum_{i=1}^{3} |\nabla b_{i}\cdot \nabla a|^{2}} > |\nabla a|^{2},
\end{equation*}
the BCD does not successfully push the spectrum downwards, even when
$\alpha$ and $\Delta E$ are small. This is exactly what happened in the previous
example: $\nabla a = (\gamma-\delta)/2, \nabla b_{3} = (\gamma+\delta)/2$, and the failure condition
reads $|(\gamma-\delta)(\gamma+\delta)| > (\gamma-\delta)^{2}$, which is always satisfied with
$\gamma, \delta > 0$.

As discussed before, because of the von Neumann theorem, non-van Hove
eigenvalue crossings must happen because of symmetries. Consider for
simplicity a two-band, 2D Hamiltonian, with a mirror symmetry
$k_{y} \mapsto -k_{y}$. This symmetry is implemented by a unitary
matrix
\begin{align*}
	H(k_{x}, -k_{y}) = U H(k_{x}, k_{y}) U^{*}.
\end{align*}
$U$ satisfies $U^{2}=1$ and therefore has eigenvalues $1$ or $-1$. The
only nontrivial case is when the eigenvalues are $-1$ and $1$. Up to a
change of basis, we can assume that $U = \sigma_{3}$, and it follows
from the Pauli matrices anticommutation relations that $b_{1}, b_{2}$
are odd with respect to the transformation $k_{y} \mapsto -k_{y}$.
Therefore, on the high-symmetry line $H(\cdot, 0)$, $b_{1}=b_{2}=0$,
and crossings happen on the codimension-1 condition $b_{3} = 0$, and
therefore generically on isolated points. At these points, the
situation is similar to the 1D case discussed in~\autoref{section:crossings},
and we expect the BCD to potentially fail
at these energies.

In 3D, the situation is even worse, because we expect crossings to
happen on curves in high-symmetry planes, and so the BCD can fail for
an extended range of energies. In practice, this range is likely to be
relatively small, as the failure of the BCD requires
$ \sqrt{\sum_{i=1}^{3} |\nabla b_{i}\cdot \nabla a|^{2}} > |\nabla a|^{2}$,
which implies that $\nabla a$ is not too large and therefore that the
crossing bands do not span a very large energy range.

We now investigate systems with symmetries to check this
condition in practice.
\subsection{Eight-bands model of graphene}
We now study the density of states (DOS) of graphene (a
three-dimensional system with two-dimensional periodicity) using an
$8\times 8$ tight-binding Hamiltonian corresponding to the $2s$ and
$2p$ orbitals (two atoms per unit cell, one $s$ orbital, three $p$
orbitals). This Hamiltonian is obtained from ab initio calculations
followed by Wannierization on a $9\times 9$ Monkhorst-Pack grid
(from Wannier90 tutorials~\cite{qiao2023projectability}). As a
reference, we use a converged DOS computed with the LT method on a
$3000 \times 3000$ grid.

Graphene has a symmetry of reflection with respect to its plane, which
leaves the two- dimensional pseudomomentum $\mathbf k$ invariant.
Therefore, for all $\mathbf k$, $H_{\mathbf k}$ is left invariant by
the mirror symmetry (obtained by flipping the signs of the $p_{z}$
orbitals). As a result, the Hamiltonian is block-diagonal, with
the $6$ $\sigma$ bands (originating from the mirror-even $s$, $p_x$,
and $p_y$ orbitals) forming the even sector and the $2$ $\pi$ bands
(derived from the mirror-odd $p_z$ orbitals) forming the odd sector.
This symmetry is even ``worse'' than the one described in the previous
section: it is a global symmetry, affecting all $\mathbf k$ points in
the same way rather than just a high-symmetry subset.

\begin{figure}[ht!]
	\centering
	\includegraphics[scale=0.29]{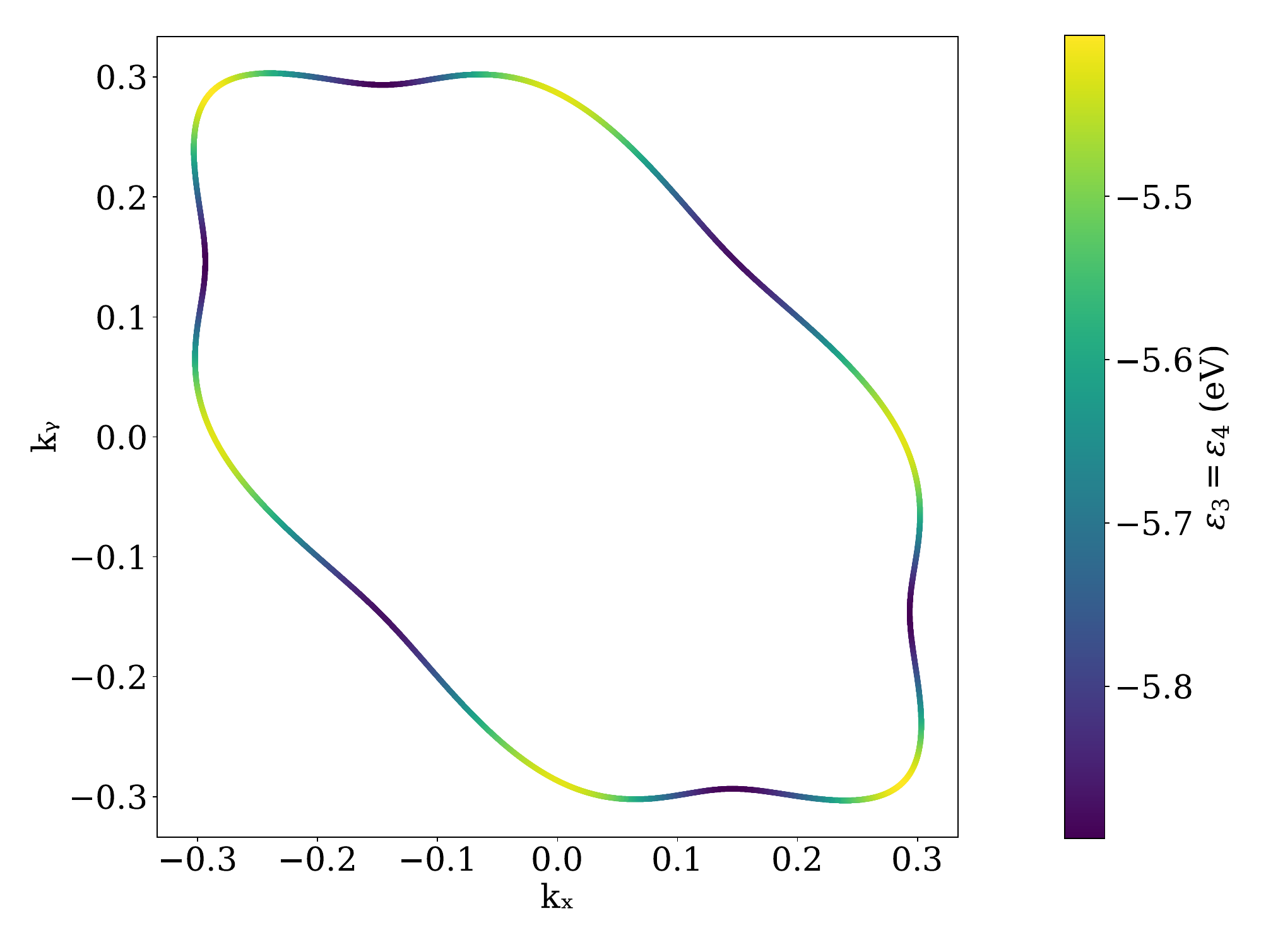}
	\caption{Crossings between the third and fourth band, protected by the mirror symmetry. The color indicates the energy of the crossing along the nodal line.}
	\label{fig:nodal_line}
\end{figure}

Since $H_{\mathbf k}$ is block-diagonal, crossings between bands from
different sectors are generically of codimension $1$, so they occur as
curves in the Brillouin zone. One such curve is illustrated
in~\autoref{fig:nodal_line}, with a crossing energy varying
continuously along the line from approximately $-5.8$ to $-5.5$ eV
(referenced from the Fermi level).

\begin{figure}[ht!]
	\centering
	\includegraphics[scale=0.29]{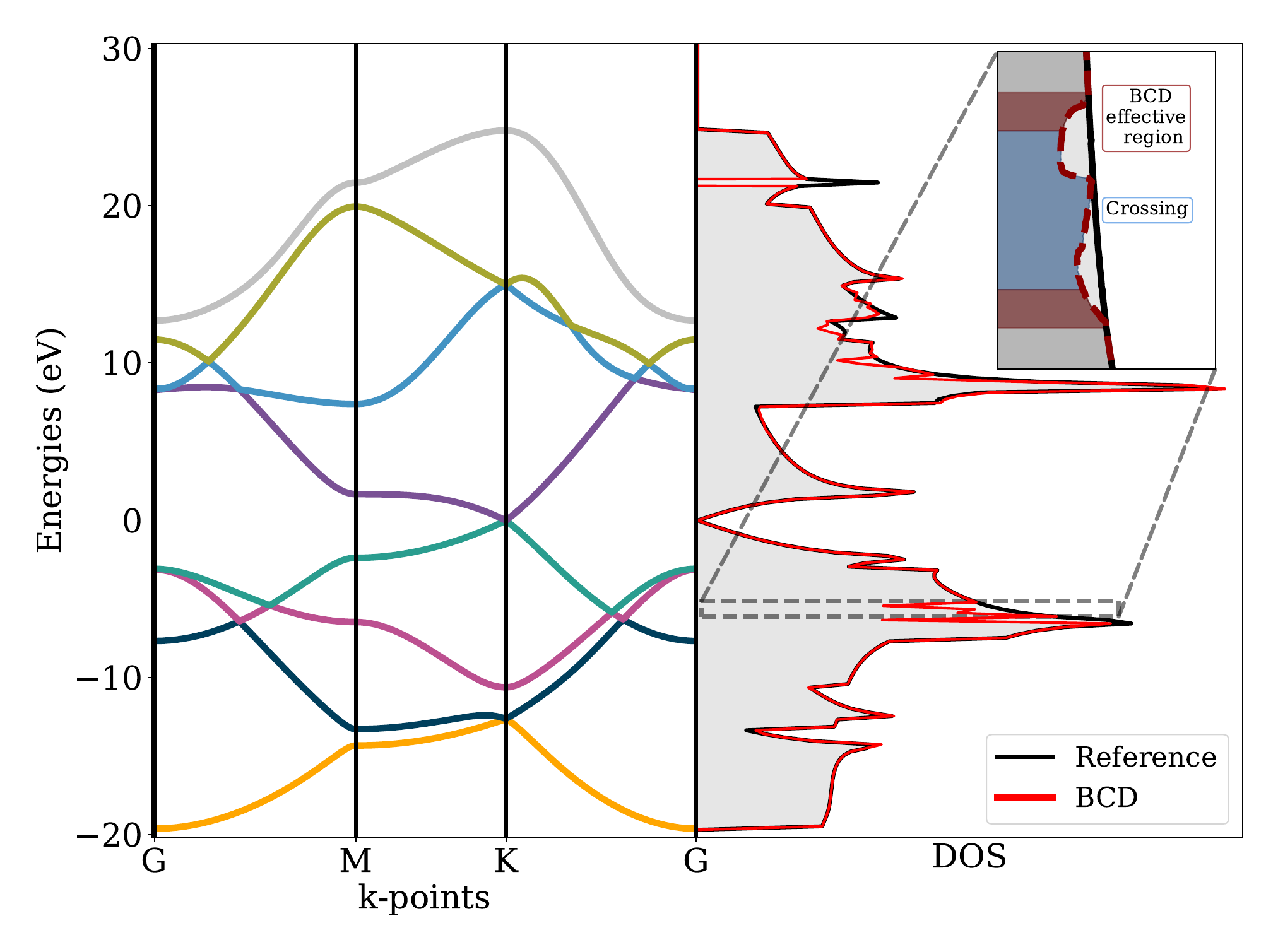}
	\caption{Bands and DOS for the 8-bands graphene with a zoom between -7.5 and 9 eV.\\ \emph{Zoom:} The blue area indicates the energy range on which the crossings from~\autoref{fig:nodal_line} occur, the dark-red range corresponds to the BCD energy smearing range for $\Delta E=0.4$ eV.}
	\label{fig:bands_dos_fullgraphene}
\end{figure}

At this type of symmetry-allowed crossings, the BCD effectively mixes
band gradients from unrelated bands, arriving at a deformation
$\mathbf h(\mathbf k)$ which has no reason to be effective for both
bands. The computation of the DOS from the BCD method in this case is
shown in~\autoref{fig:bands_dos_fullgraphene}. As in the 1D case, in
the energy range spanned by the crossing curves, the BCD method
systematically underestimates the DOS, with the discrepancy
concentrated in bands of width $\sim\Delta E$ around each crossing
energy. Other crossings higher in energy produce the same effect.

\subsection{The free electron gas (free Laplacian)}
Consider the free electron gas, obtained by artificially considering the completely homogeneous operator $-\Delta$ as a periodic operator with lattice $2\pi
	\mathbb{Z}^d$. In a sense, this is the most symmetric model possible, having all possible geometric symmetries. Its eigenvalue bands can be labelled by the index $\mathbf
	G \in \mathbb{Z}^d$ as $\varepsilon_{\mathbf G \mathbf k} =
	|\mathbf{k}+\mathbf{G}|^{2}$ and can be seen as shifted paraboloids. This model has a single van Hove singularity, at the spectrum edge $0$.

\subsubsection{1D electron gas}
In one dimension, the band structure of the free Laplacian consists of parabolic branches
\begin{equation*}
	\varepsilon_G(k) = (k+G)^2, G\in\mathbb Z, k\in\left[-\tfrac12,\tfrac12\right].
\end{equation*}
Two types of crossings occur: at $k=0$ with energies $G^2$ ($G\neq0$),
and at $k=\pm 1/2$ with energies $(1/2+G)^2$ ($G\in\mathbb Z$). The band
structure exhibits reflection symmetries about $k=0$ (since
$\varepsilon_{-G}(-k)=\varepsilon_G(k)$) and about $k=1/2$ (since
$\varepsilon_{1-G}(1/2-k)=\varepsilon_G(1/2+k)$). Consequently, at these symmetric points
the deformation field $h(k)$ vanishes identically for any choice of
$E$ and $\Delta E$; the BCD therefore does nothing at these isolated $k$
points. This is the same case as the two independent bands
(c.f~\autoref{section:crossings}) but with equal and opposite
gradients i.e $\gamma = \delta$. Therefore we expect the convergence of the density of states with
respect to the number of quadrature points $N$ to be slow at the
crossing energies $(1/2)\mathbb Z$. Away from these special energies,
a single band dominates and the BCD is effective provided $\Delta E$ is
chosen small enough. However, for a given finite $\Delta E$, the BCD can
still fail for some energies close to the crossing energies.

\begin{figure}[hbtp]
	\centering
	\subfloat{\includegraphics[scale=0.24]{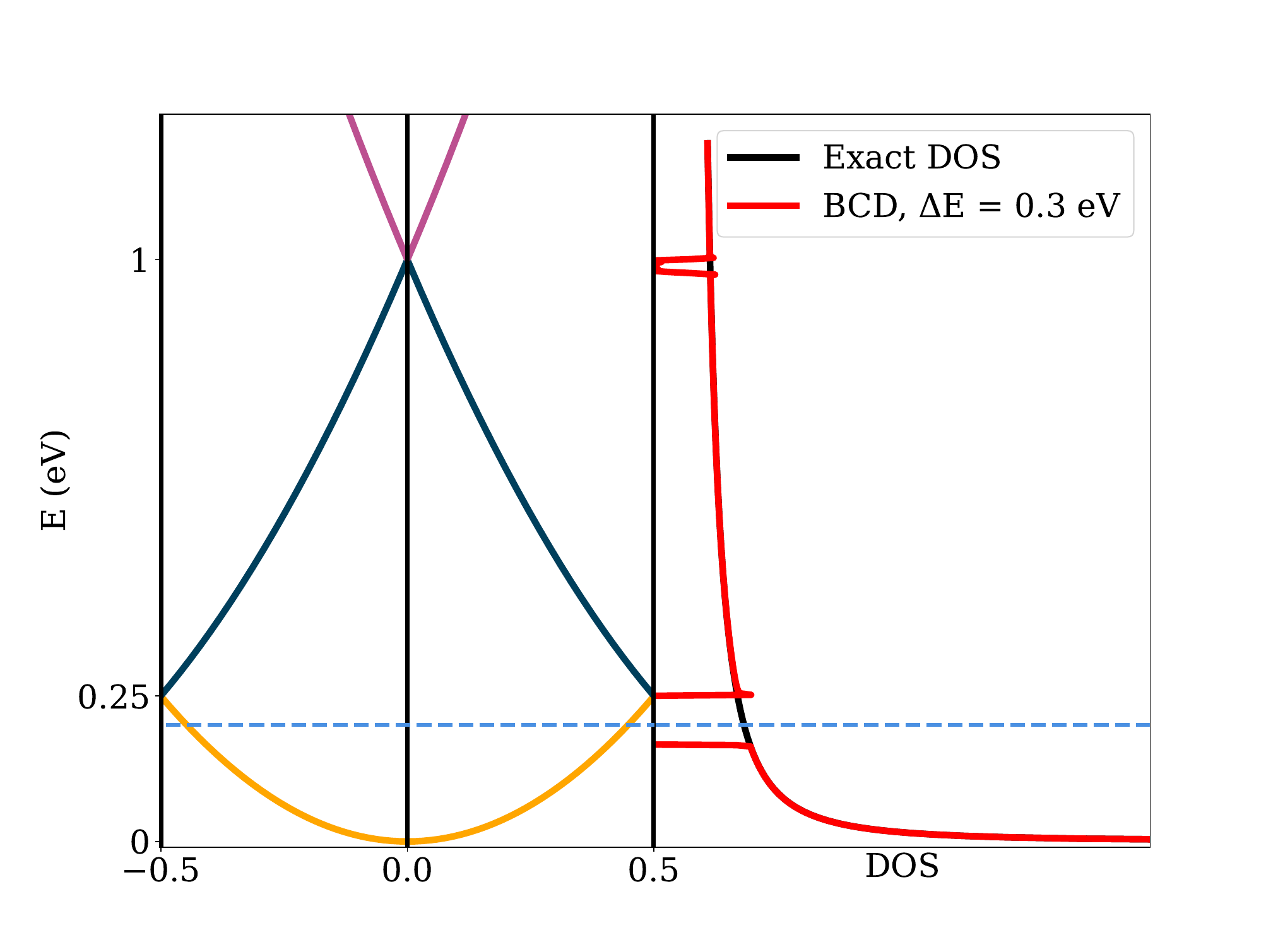}}
	\subfloat{\includegraphics[scale=0.24]{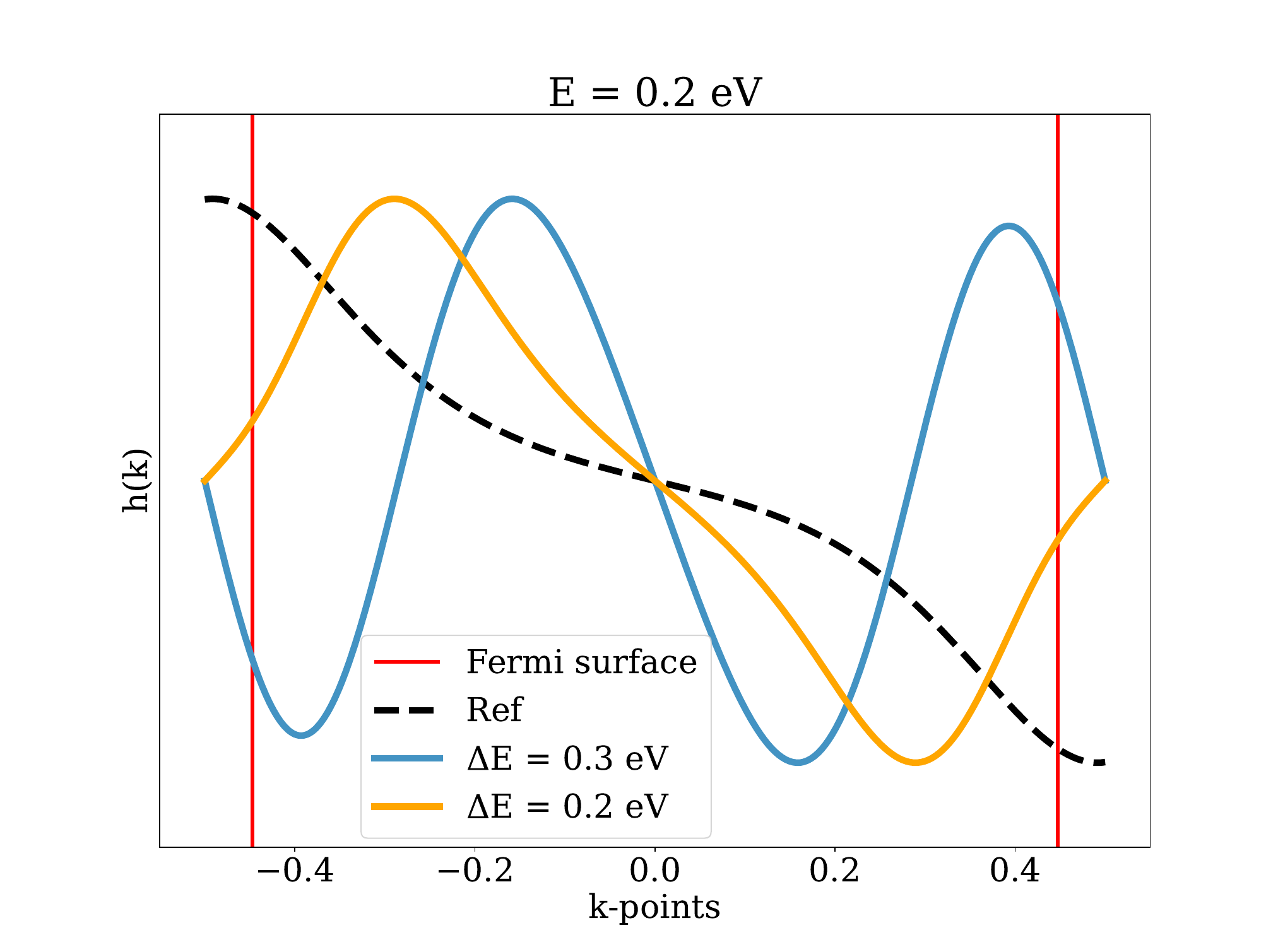}}
	\caption{\emph{Left panel:} Three first bands ($G_1=0$, $G_2=-1$, $G_3=1$ ) with corresponding exact DOS $E^{-1/2}$ and DOS computed with the BCD for $\alpha=0.1$ and $\Delta E= 0.3$ eV. The blue line indicates $E=0.2$ eV. \\
		\emph{Right panel:} Deformation field for the 1D free electron gas at $E=0.2$ eV with parameters $\Delta E=0.3$ and $\Delta E=0.2$ eV. The reference deformation field is computed by taking into account only the first band with $\Delta E=0.2$ eV.}
	\label{fig:1D_gas}
\end{figure}

We illustrate this process in~\autoref{fig:1D_gas}. The left panel
shows the DOS computed by the BCD method: it fails in energies smaller
than each crossing energy (e.g., $0.25$ eV and $1$ eV). The right panel displays the deformation field $h(k)$ as a function of $k$ for several values of $\Delta E$. For large $\Delta E$, the second band introduces a spurious contribution that leads to a deformation pushing the eigenvalues upwards instead of downwards. As $\Delta E$ decreases, the BCD successfully pushes the spectrum downward. We also show a reference deformation field computed by including only the first band, which is the only band crossing the Fermi surface at this energy. For $\Delta E = 0.3\,\mathrm{eV}$, the second band contribution reverses the sign, leading to an upward deformation. For $\Delta E = 0.2\,\mathrm{eV}$, its contribution is smaller than that of the first band, and the BCD yields a downward deformation.

To be more precise, consider two bands crossing with indices
$G' \neq G''$. The crossing is located at $k_{0}= -(G' + G'')/2$. We
introduce the shifted variable $k=q+k_{0}$, to bring the crossing point
to $q=0$. Setting $G= G'-G'' \in \mathbb Z$, the two bands take the form
\begin{equation*}
	\varepsilon_{\pm}(q)=\left(q\pm\frac G2\right)^2,
\end{equation*}
whose common energy value at the crossing is $E_0=(G/2)^2$ eV. For
$E < E_{0}$, the band at the Fermi level is $\varepsilon_{-}$, and its derivative
goes from positive for $q < 0$ to negative for $q > 0$. The correct
deformation $h$ should then be negative for $q < 0$ and positive for
$q > 0$.

For $E$
close to $E_0$, so that the effect of the other bands can be ignored, the deformation field reads
\begin{equation*}
	h(q)=\varepsilon_{-}'(q)\,\exp\left(-\left(\frac{\varepsilon_{-}(q)-E}{\Delta E}\right)^2\right)
	+\varepsilon_{+}'(q)\,\exp\left(-\left(\frac{\varepsilon_{+}(q)-E}{\Delta E}\right)^2\right).
\end{equation*}
By symmetry, $h(0)=0$. Differentiating with respect to $q$ yields
\begin{equation*}
	h'(0)= 4\exp\left(-\frac{(E-E_0)^2}{\Delta E^2}\right)
	\left(1+\frac{G^2(E-E_0)}{\Delta E^2}\right).
\end{equation*}
Fix $E$ close to $E_0$ but smaller. For $\Delta E$ sufficiently small, the
BCD is successful and $h'(0) > 0$. However, when $\Delta{E}$ is
increased so that $E_0-(\Delta E/G)^2<E<E_0$, $h'(0)$ becomes negative; the
deformation pushes in the wrong direction and the BCD fails.

\subsubsection{2D and 3D}
In two and three dimensions, the Fermi surface exhibits a richer structure.
For a generic energy $E$, the sheets
$S_{\mathbf{G}}(E) = \{\mathbf{k},\; |\mathbf{k}+\mathbf{G}|^{2} = E\}$
are respectively circles and spheres.
For $E<0.25$ eV there are no crossings because only the band $\varepsilon_{\mathbf{0}}$ lies in the Brillouin zone; hence the BCD performs well in this regime up to a $\sim \Delta E$ neighbourhood before $E=0.25$ eV.
Starting at this energy, crossings appear generically: first pairwise, then with three or more circles.

\begin{figure}[ht!]
	\centering
	\subfloat{\includegraphics[scale=0.46]{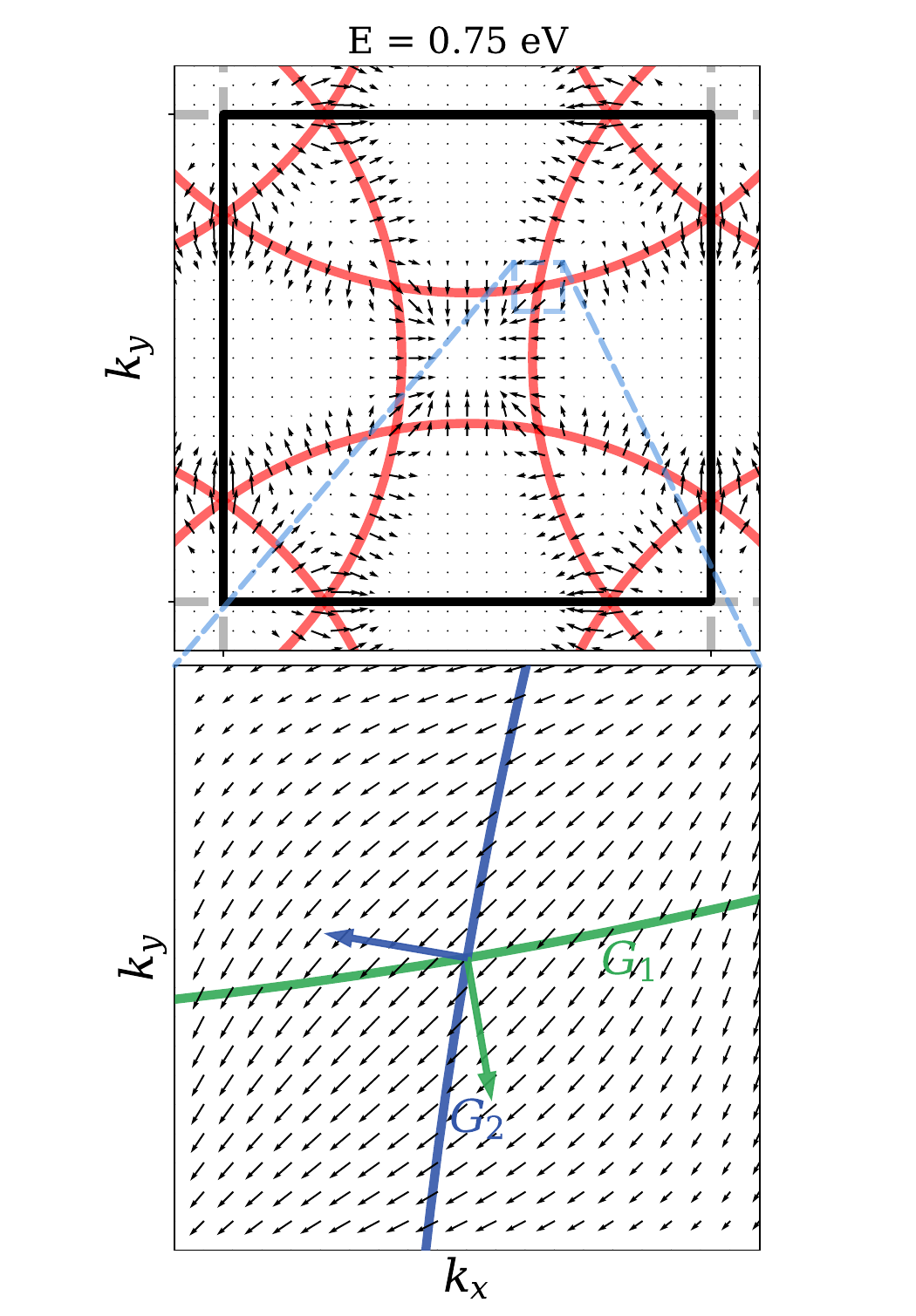}}
	\subfloat{\includegraphics[scale=0.46]{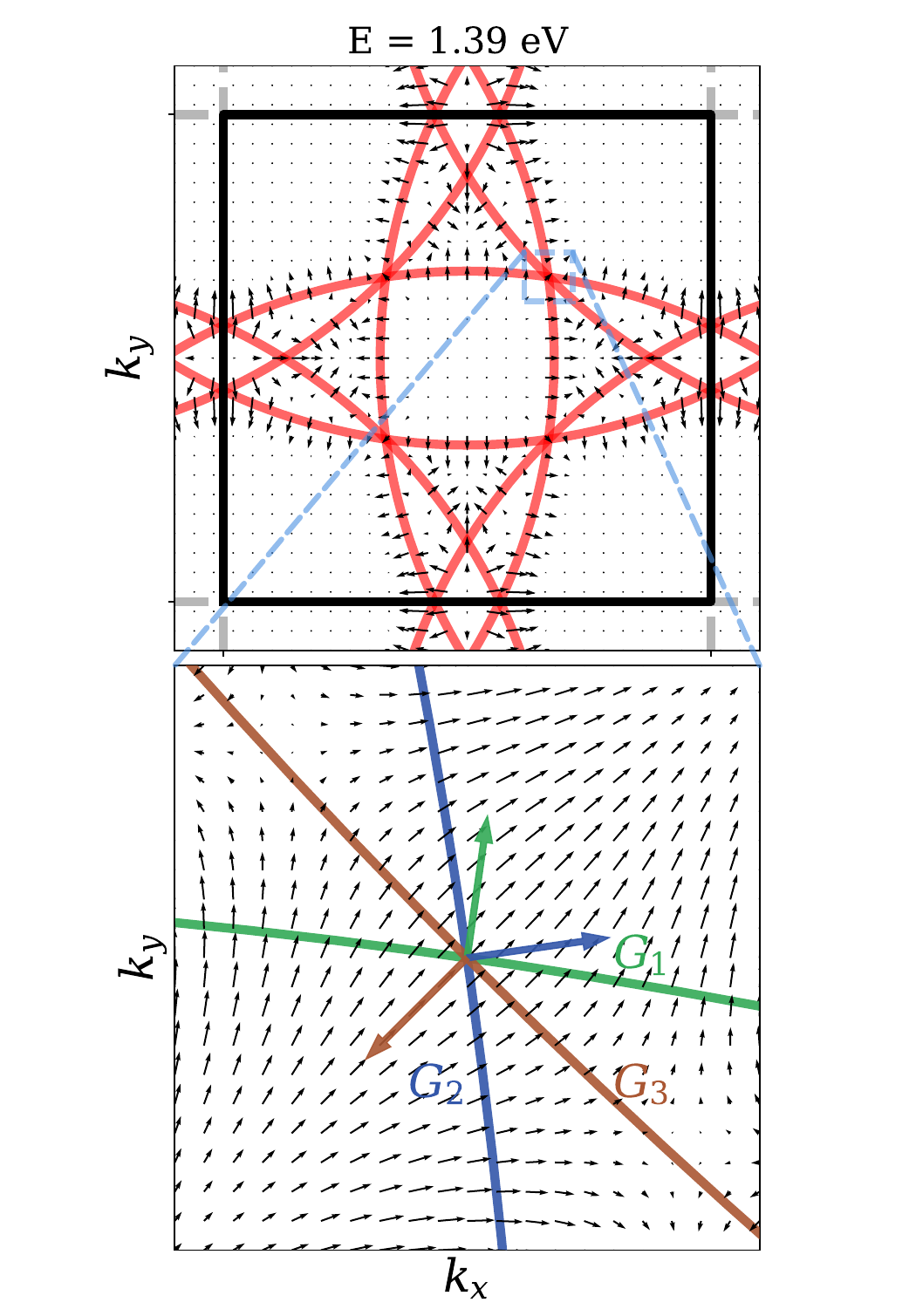}}\hfill
	\caption{\emph{Left panel:} Fermi surface for $E=0.75$ eV. At this
		exact energy the crossings cause no underestimation of the DOS thanks to symmetries.\\
		\emph{Right panel:} Fermi surface for $E=25/18 \approx 1.39$ eV. Triple crossings appear at the points $\mathbf{k}=(\pm 1/6,\pm 1/6)$ with associated lattice points $\mathbf{G}_1=\mp(0,1)$, $\mathbf{G}_2=\mp(1,0)$ and $\mathbf{G}_3=(\pm 1,\pm 1)$. Here the deformation pushes the spectrum in the wrong direction on $S_{\mathbf{G}_3}$.}
	\label{fig:2D_gas_crossings}
\end{figure}

\begin{figure}[ht!]
	\centering
	\includegraphics[scale=0.29]{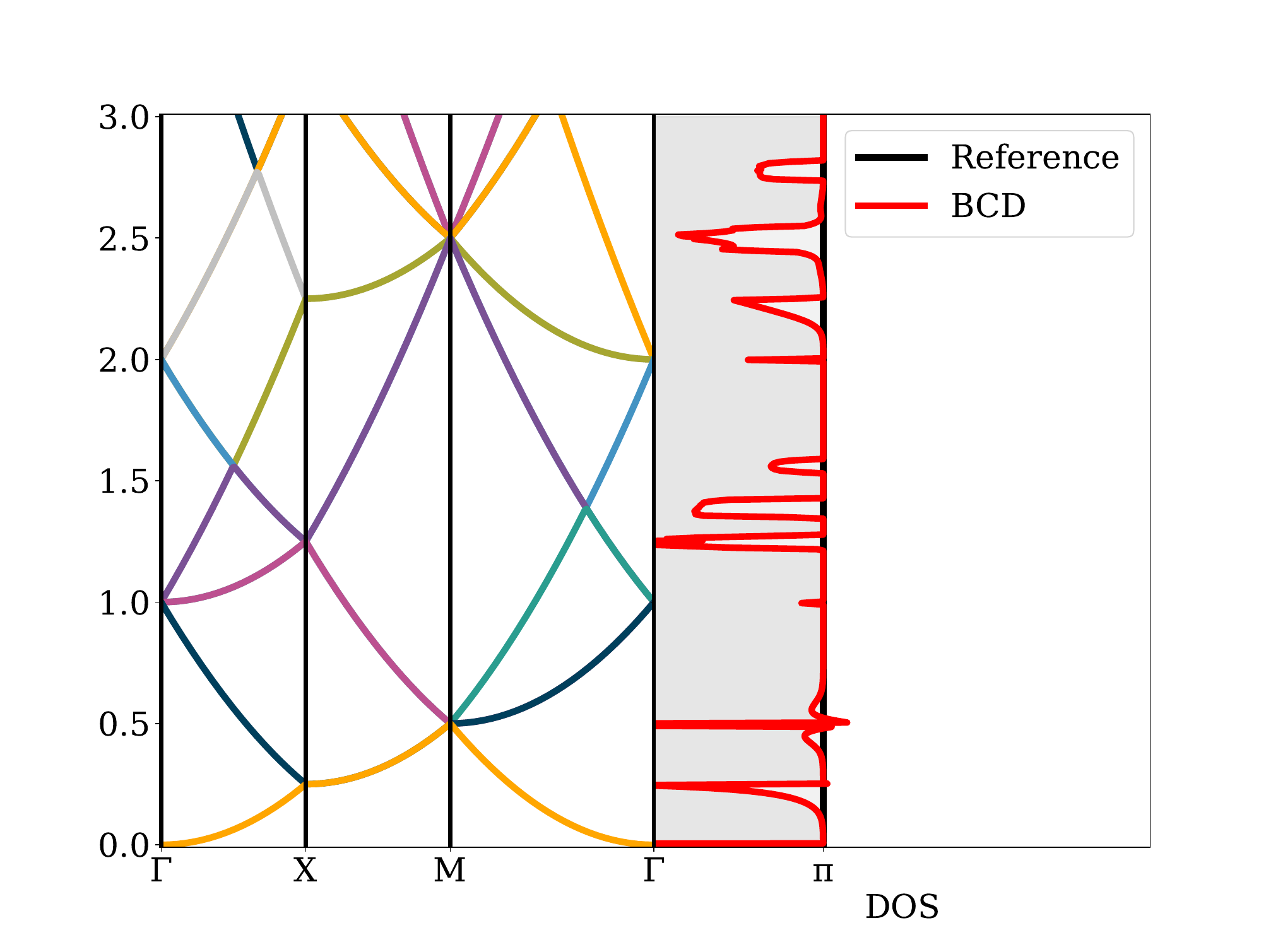}
	\caption{Bands and DOS of the 2D electron gas computed using the BCD with parameters $\alpha=0.1$ and $\Delta E=0.1$ eV. }
	\label{fig:2D_gas_dos}
\end{figure}

\paragraph{Crossings of two circles.}
Consider a crossing at $\mathbf{k}$ between $S_{\mathbf{G}}$ and $S_{\mathbf{G}'}$.
For sufficiently small $\Delta E$, the BCD deformation field is well approximated by the average of the two normals:
\begin{equation*}
	\mathbf{h}(\mathbf{k}) \approx -\alpha\bigl(\nabla\varepsilon_{\mathbf{G},\mathbf{k}}+\nabla\varepsilon_{\mathbf{G}',\mathbf{k}}\bigr) = -2\alpha\bigl((\mathbf{k}+\mathbf{G})+(\mathbf{k}+\mathbf{G}')\bigr).
\end{equation*}
To first order in $\alpha$, the complex shift of the band $\varepsilon_{\mathbf{G}}$ is determined by
\begin{equation*}
	\nabla\varepsilon_{\mathbf{G},\mathbf{k}}\cdot\mathbf{h}(\mathbf{k}) = -2\alpha\bigl(|\mathbf{k}+\mathbf{G}|^{2}+(\mathbf{k}+\mathbf{G}')\cdot(\mathbf{k}+\mathbf{G})\bigr).
\end{equation*}
Since $|\mathbf{k}+\mathbf{G}|=|\mathbf{k}+\mathbf{G}'|$, the
Cauchy–Schwarz inequality implies that this quantity is negative
unless $(\mathbf{k}+\mathbf{G}) = -(\mathbf{k}+\mathbf{G}')$. This
condition holds precisely when $2\mathbf{k}=-(\mathbf{G}+\mathbf{G}')$,
i.e. when $\mathbf{k}\in (1/2)\mathbb{Z}^{d}$. The same reasoning
applies to the band $\varepsilon_{\mathbf{G}'}$. In contrast with the
eight-band graphene model, the gradients have equal magnitude;
consequently, the deformation either acts in the correct direction or
vanishes, but does not shift the eigenvalues in an incorrect
direction. This is illustrated in the left-panel
of~\autoref{fig:2D_gas_crossings}. However, for $E$ in a $\sim \Delta E$ neighborhood of the crossing energy, the BCD is incorrect, as can be
seen in \autoref{fig:2D_gas_dos}.

\paragraph{Crossings of three circles.}
A similar reasoning can be applied to crossings of $3$ circles with
centers $-\mathbf G_{i}, i=1, \dots, 3$. In this case,
\begin{equation*}
	\nabla\varepsilon_{\mathbf{G}_1,\mathbf{k}}\cdot\mathbf{h}(\mathbf{k}) = -2\alpha\left(|\mathbf{k}+\mathbf{G}_1|^{2}+\sum_{i=2 }^3(\mathbf{k}+\mathbf{G}_1)\cdot(\mathbf{k}+\mathbf{G}_i)\right).
\end{equation*}
Using Cauchy–Schwarz, the absolute value is bounded by
$2\times 3\alpha E$, but the sign can be either positive or negative. In
particular, when the number of bands is odd, the sum
$\sum_i(\mathbf{k}+\mathbf{G}_0)\cdot(\mathbf{k}+\mathbf{G}_i)$ can be
sufficiently negative to flip the sign of the imaginary shift, pushing
the spectrum in the wrong direction.
We give an example of such crossing and its effect on the deformation
in the right panel of~\autoref{fig:2D_gas_crossings}.

Such crossings between $n$ circles are non-generic for $n > 2$: they
happen at isolated energies. However, for all choices of
$\mathbf G_{1}, \mathbf G_{2}, \mathbf G_{3}$ which are not aligned, one
can find an energy at which the relevants bands cross (the square of
the distance from all $\mathbf G_{i}$ to the circumcenter of the
triangle). Therefore, there is an infinite number of failure points of
the BCD. Furthermore, these crossings get denser and denser as one
increases the energy, and the BCD gives consistently wrong
results, as can be seen in~\autoref{fig:2D_gas_dos}.

\begin{figure}[ht!]
	\centering
	\includegraphics[scale=0.29]{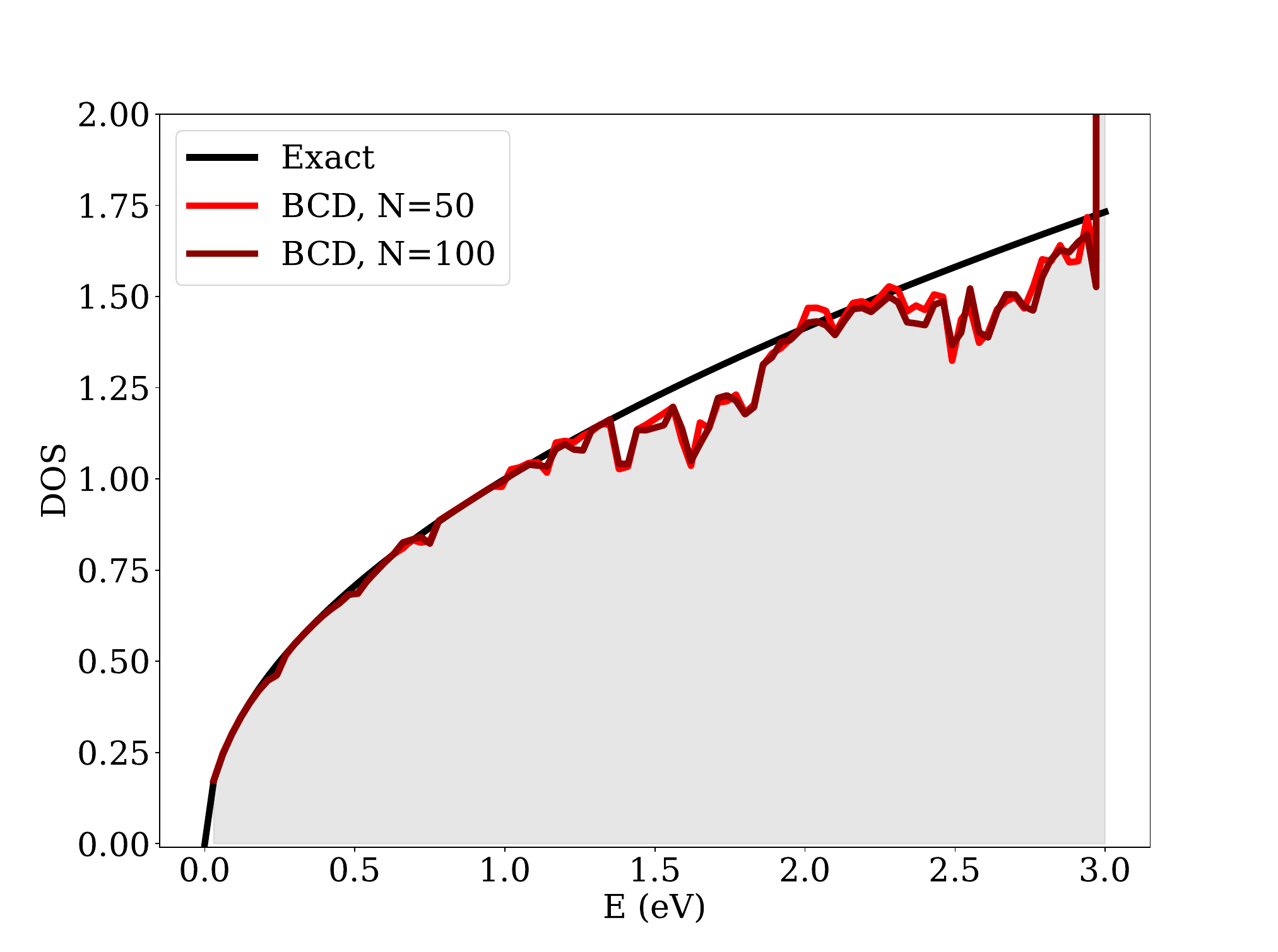}
	\caption{DOS of the 3D electron gas computed using the BCD with parameters $\alpha= 0.1$ and $\Delta E=0.1$ eV for $N=50$ and $N=100$.}
	\label{fig:3D_gas_dos}
\end{figure}

\paragraph{3D free electron gas.}
In three dimensions, the crossings become even more numerous, and the
BCD is unusable. The DOS present failures for $E> 0.25$ eV resulting in an underestimation of the DOS, with complete failure for all energies $E> 1$ eV as illustrated in~\autoref{fig:3D_gas_dos}. However, there are still errors attributable to the discretization method, as illustrated by considering both $N=50$ and $N=100$.

\subsection{Discussion}
One can wonder if the failure of the BCD is fundamental or stems from
a wrong choice of $\mathbf h(\mathbf k)$ in the multiband case.
Indeed, the choice of such a deformation has been considered before in
a related setting, that of choosing a conjugated operator to prove
scattering properties for periodic operators in the framework of
Mourre theory \cite{gerard1998mourre}. Indeed, the choice of conjugate
operator in \cite{gerard1998mourre} can be thought of as a version of
the BCD for infinitesimally small $\alpha$ (the first order of the
expansion of $H_{\mathbf{k} + \alpha \mathbf h(\mathbf k)}$ involves a
commutator of $H$ with the Mourre conjugate operator). In
\cite{gerard1998mourre}, the authors make a particular choice of
$\mathbf h(\mathbf k)$ which they prove to work (in the sense of
satisfying the Mourre inequality, which is analogue for us to pushing
the spectrum downwards) for \emph{almost every} energy. However, their
construction relies on very delicate geometry arguments (Whitney
stratification of semianalytic sets) that are hard to turn into a
numerical method. Furthermore, even this theory must exclude a locally
finite number of energies: for instance, it is ineffective on the
Hamiltonian $
	\begin{pmatrix}
		k & 0 \\0& -k
	\end{pmatrix}
$. A more sophisticated method would attempt to either detect this
block-diagonal structure, or operate natively to deform different
components of the Hamiltonian in different $\mathbf k$-space
directions. We have been unable to find such a general procedure,
however.

All the examples of failure of the BCD in this section were artificial.
In particular, all possessed a ``block diagonal'' structure which
would make it trivial to compute the DOS of each subblock using the
BCD method and summing up the final result. Realistic models
of materials (with full periodicity, i.e. excluding surfaces which are
3D models with 2D periodicity) are not expected to have this structure
globally; indeed, this would contradict the conjectured irreducibility
of the Bloch variety of elliptic periodic operators
\cite{kuchment2023analytic}. However, from the discussion of
\autoref{sec:crossings_symmetries}, they could very well have this
block diagonal structure on high-symmetry lines or planes, where the BCD could fail. Fortunately, this phenomenon
appears somewhat rare among ``simple'' materials, and in particular is
not present in the systems we test in the rest of this paper, as well
as a number of realistic bulk 3D materials (we tested on silicon,
calcium, copper, and lead). This might be because such problematic
crossings only appear when bands are subject to strong variations
(such as band inversion),
which is not the case in simple materials.

To detect a failure of the BCD, a simple and inexpensive diagnostic is
to monitor the imaginary part of the eigenvalues of the deformed
Hamiltonian $H_{\mathbf{k}+i\mathbf{h}(\mathbf{k})}$ to ensure that
the eigenvalues are pushed down instead of up. We report failure when
\begin{equation}
	\exists n, \mathbf{k}: \quad |\varepsilon_{n\mathbf{k}} - E| < \frac{\Delta E}{2} \text{ and } \mathrm{Im}(\varepsilon_{n,\mathbf{k}+i\mathbf{h}(\mathbf{k})}) > \mathrm{tol}. \label{eq:diagnostic}
\end{equation}
The tolerance parameter accounts for the wanted accuracy and numerical
errors in the eigenvalue computation, typically set to
$\mathrm{tol} \sim 10^{-6}$ eV for high-precision DOS calculations. A positive
imaginary part signals a failure of the BCD and an underestimation of
DOS; in this case, an alternative method must be used.

In the following, having no non-artificial failure cases, we focus
on cases where the BCD works. A more systematic study would be
interesting research but is outside of the scope of this paper.

\section{Methods}
Here we discuss the numerical setting in which we did our benchmark and the main computational differences between the methods.
\paragraph{Tight-binding Hamiltonian}
We consider a tight-binding Hamiltonian of the form
\begin{equation*}
	H_\mathbf{k}\propto\sum_{\mathbf{R} \in \mathcal{R}_N} \mathrm{e}^{-i \mathbf{k}\cdot\mathbf{R}}H_{\mathbf{R0}}.
\end{equation*}
This is a truncated Fourier series over a set $\mathcal{R}_N \subset \mathcal{R}$ of given
$\mathbf{R}$ vectors and corresponding $H_\mathbf{R}$ matrices coupling
different Wannier functions from different unit cells.
This Hamiltonian is obtained either from an ansatz (e.g., the two-band
model of graphene), or from ab-initio computations like wannierization (e.g., $\mathrm{SrVO_3}$).

\paragraph{Computation time}
In the case of the PTR and the LT methods, $H_{\mathbf{k}}$ is to be
evaluated on a cubic grid, so that this sum can be performed using the
fast Fourier transform (if memory permits). In the case of the IAI
method, the grid is not fully regular, but is sufficiently tensorized
(i.e. a significant portion of points share one or two coordinates)
that some computations can be reused using partial Fourier transforms,
cutting down on the total evaluation time of the
Hamiltonian~\cite[Appendix C]{kaye2023automatic}. For the BCD method
by contrast, the (deformed) grid is fully general with no particular
structure. Furthermore, the computations need to be performed in
complex arithmetic, and the computation of several derivatives is
required. In our benchmarks, utilizing the optimized implementations
described in~\cite{kaye2023automatic}, we observe that for challenging
three-dimensional applications, for a single fixed energy $E$ and a
fixed total number of evaluation points, the BCD method is the most
computationally intensive. It is approximately 80-100 times slower
than PTR and IAI, mainly due to the lack of structure in the deformed
grid, which prevents the use of fast Fourier transforms and requires
more expensive interpolation methods. The difference is less
significant in 1D and 2D, where tensorization effects are weaker. Lastly, for
reference, we compare with a barebones implementation of the LT method
which is not optimized and performs 3-5 times slower than PTR and IAI
(for a given number of Hamiltonian evaluations).

\begin{figure}[ht!]
	\centering
	\includegraphics[scale=0.29]{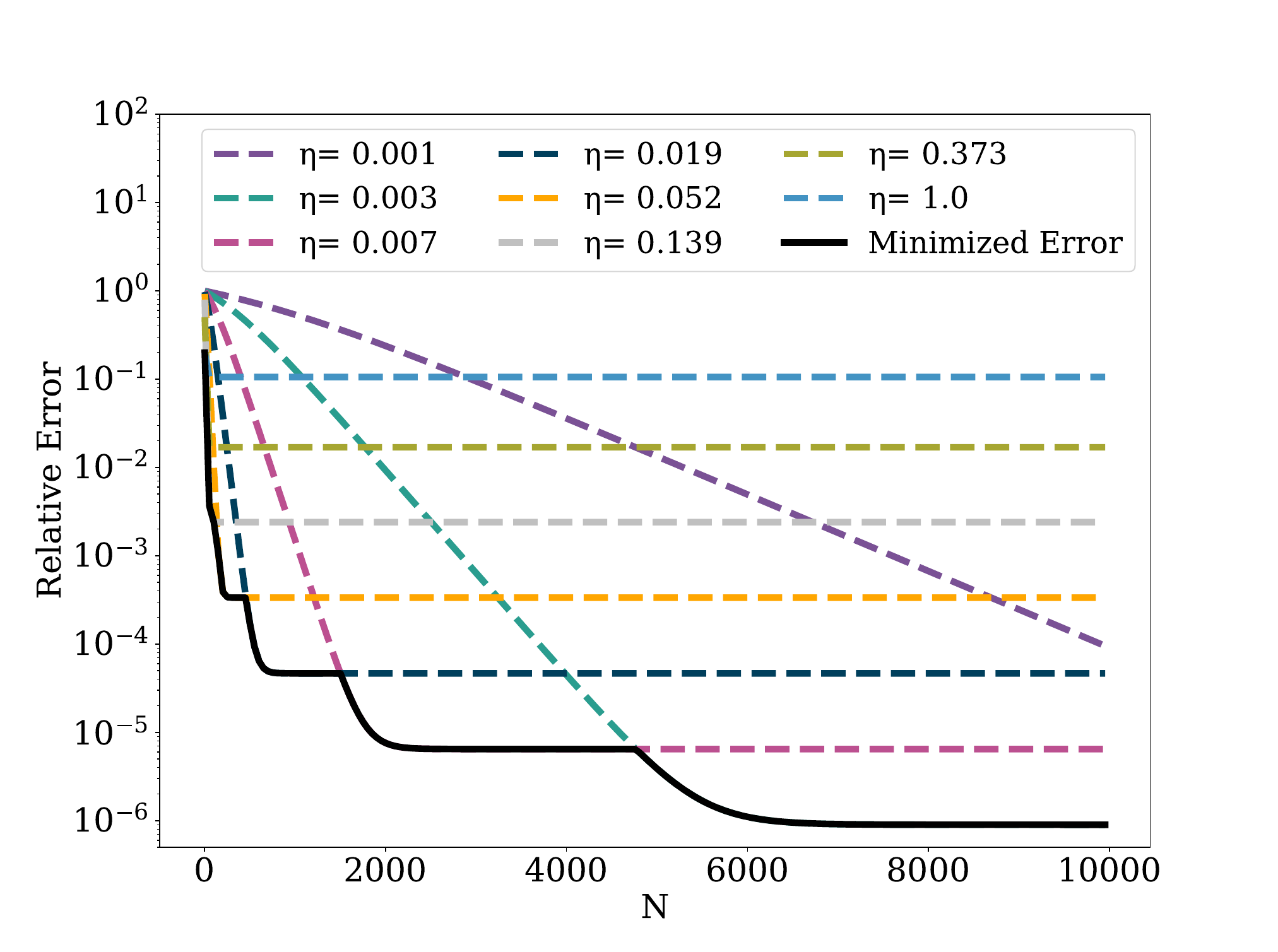}
	\caption{Choosing $\eta$ to minimize the error to the
		reference zero-smearing DOS for the 1D monatomic chain, using PTR.}
	\label{fig:multipleeta}
\end{figure}

\paragraph{Choice of smearing}
For the IAI and PTR methods, we have to select the smearing parameter
$\eta$. For a fixed number of quadrature points, when $\eta$ is
reduced, the regularization error goes down, but the discretization
error goes up. The optimal $\eta$ must therefore go to zero as the number
of quadrature points is increased. To ensure a best-case scenario for the smearing
methods, in our numerical tests, for a fixed computational budget
(number of quadrature points), we selected an almost-optimal $\eta$ by
computing the result with different $\eta$, and choosing the one that achieved
the minimum error, as measured by a comparison with a converged
reference result. We illustrate this in~\autoref{fig:multipleeta}.
\section{Smearing-free DOS results}\label{section:results}
\subsection{Monatomic chain}
\begin{figure}[ht!]
	\centering
	\includegraphics[scale=0.29]{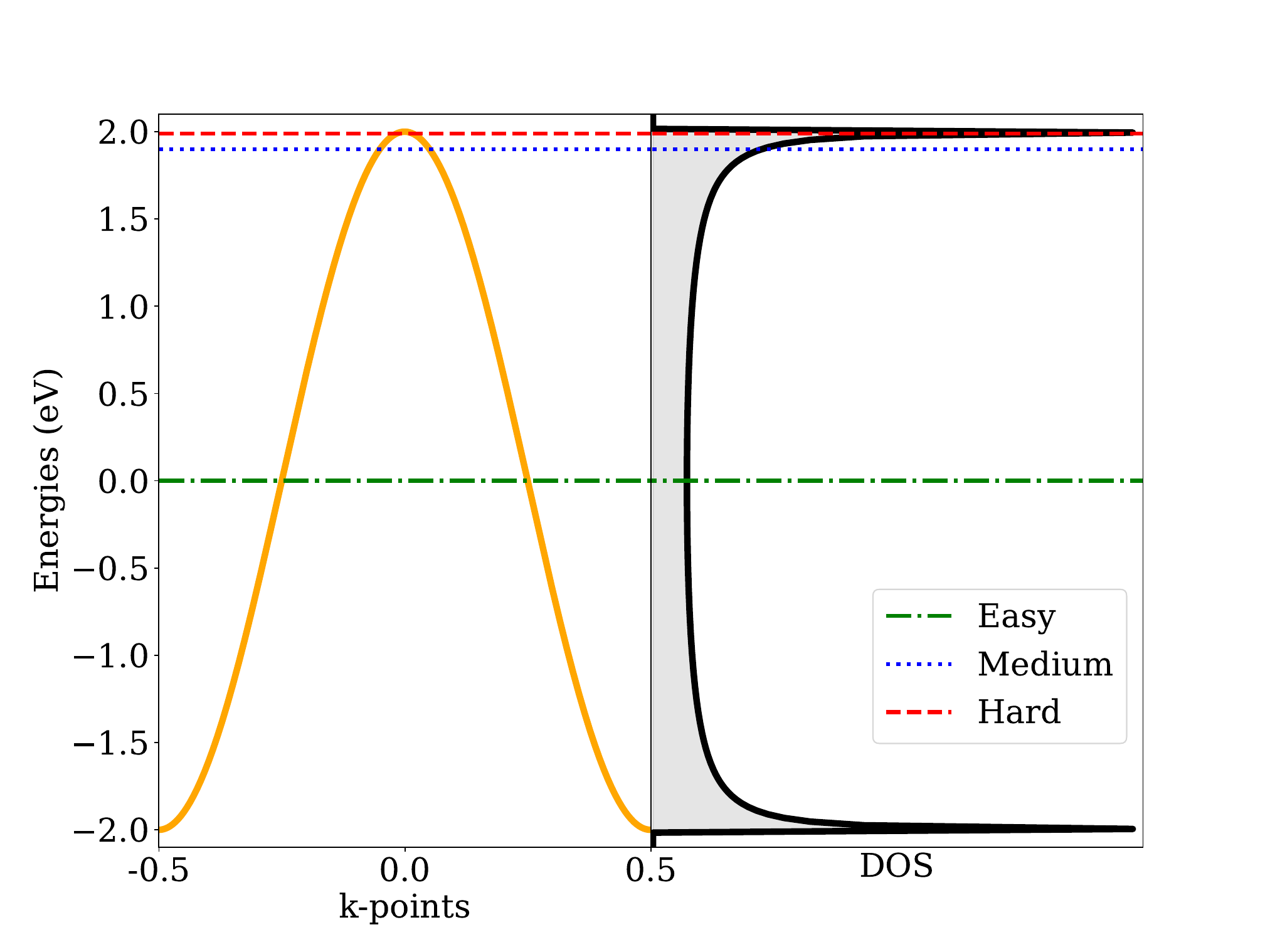}
	\caption{Bands and DOS of the 1D monatomic chain with hopping constant
		$t=1$.}
	\label{fig:dosmono1d}
\end{figure}

We consider a tight-binding model for a monatomic chain with hopping coefficient $t$. This model serves as a benchmark test case since both its Hamiltonian and DOS have explicit analytical form. The real-space Hamiltonian $H$, acting on $\ell^2(\mathbb{Z})$, is the infinite matrix
\begin{equation}
	H=\begin{pmatrix}
		\ddots & \ddots & \ddots &        &        &        \\
		       & t      & 0      & t      &        &        \\
		       &        & t      & 0      & t      &        \\
		       &        &        & \ddots & \ddots & \ddots
	\end{pmatrix}.
\end{equation}
Applying the Bloch transform to $H$ results in
\begin{equation}
	H_{k}= 2t\cos(k),
\end{equation}
with the Brillouin zone $\mathcal{B} = 2\pi\mathbb{R}  / \mathbb{Z}$.
The exact DOS is derived from~\autoref{eq:dosgrad}
\begin{equation}
	\left\{\begin{matrix}
		\displaystyle{D^{\rm 1D}_{\mathrm{ref}}(E)=\frac{1}{\pi\sqrt{(2t)^2-E^2}}}, & E \in (-2t,2t),     \\
		0,                                                                          & \mathrm{elsewhere}.
	\end{matrix}\right.
	\label{eq:dos1D}
\end{equation}
The DOS exhibits two van Hove singularities at the band edges located
at energies $E_\pm=\pm 2 t$. As is typical for one-dimensional
system~\cite{vanhove1953occurrence}, near the van Hove singularities the DOS behaves as
$(E_\pm \mp E)^{-1/2}$. The bands and the corresponding DOS are shown
in~\autoref{fig:dosmono1d}.

We benchmarked the numerical methods for three representative energies
depending on their relative distance to the van Hove
singularities: an easy case, in a smooth region ($E=0$ eV), a medium case,
somewhat close to a van Hove singularity ($E=1.9$ eV), and a hard case,
very close to a van Hove singularity ($E=1.99$ eV).

\begin{figure}[ht!]
	\centering
	\subfloat[Easy]{\includegraphics[scale=0.24]{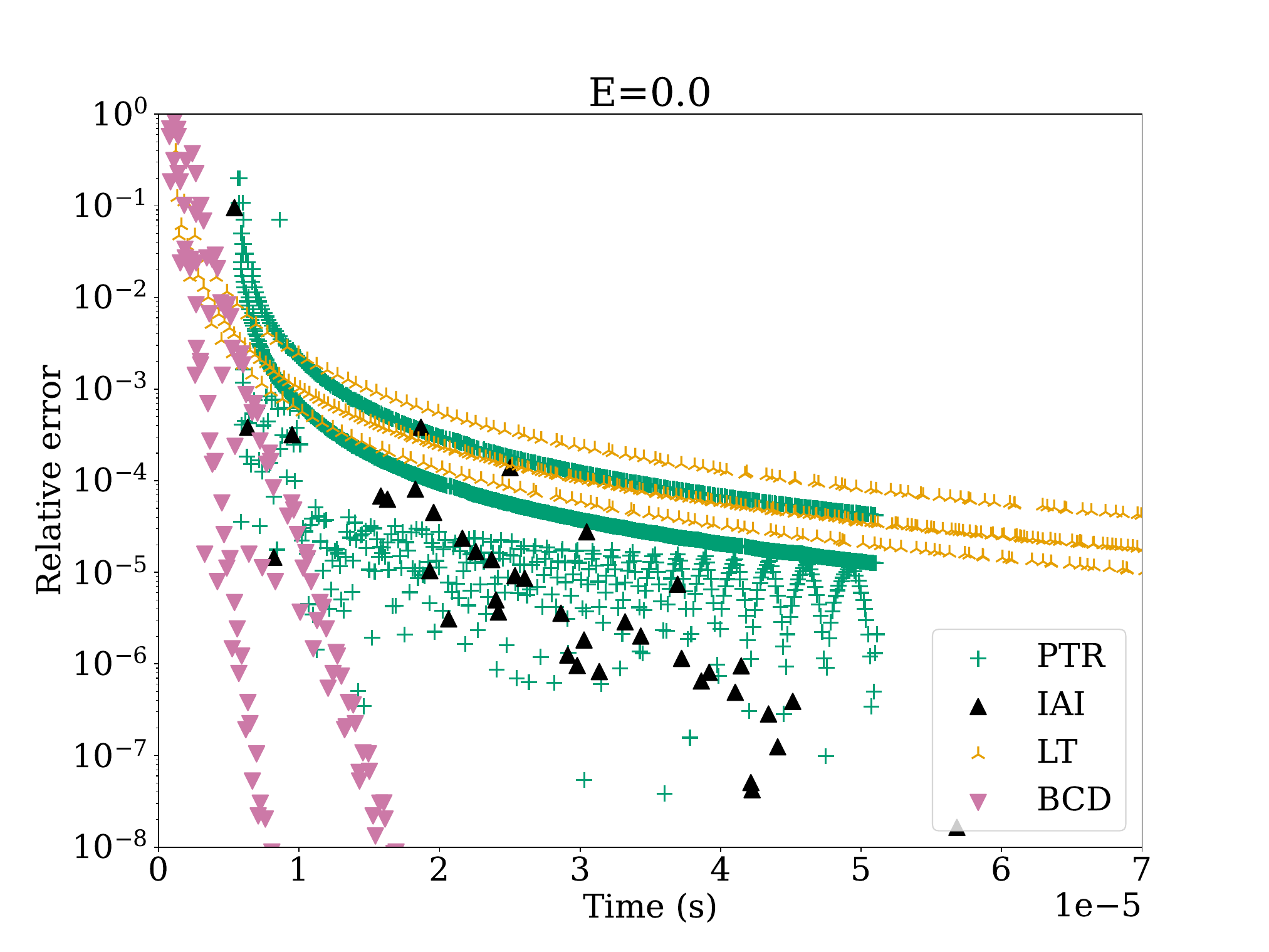}}
	\subfloat[Medium]{\includegraphics[scale=0.24]{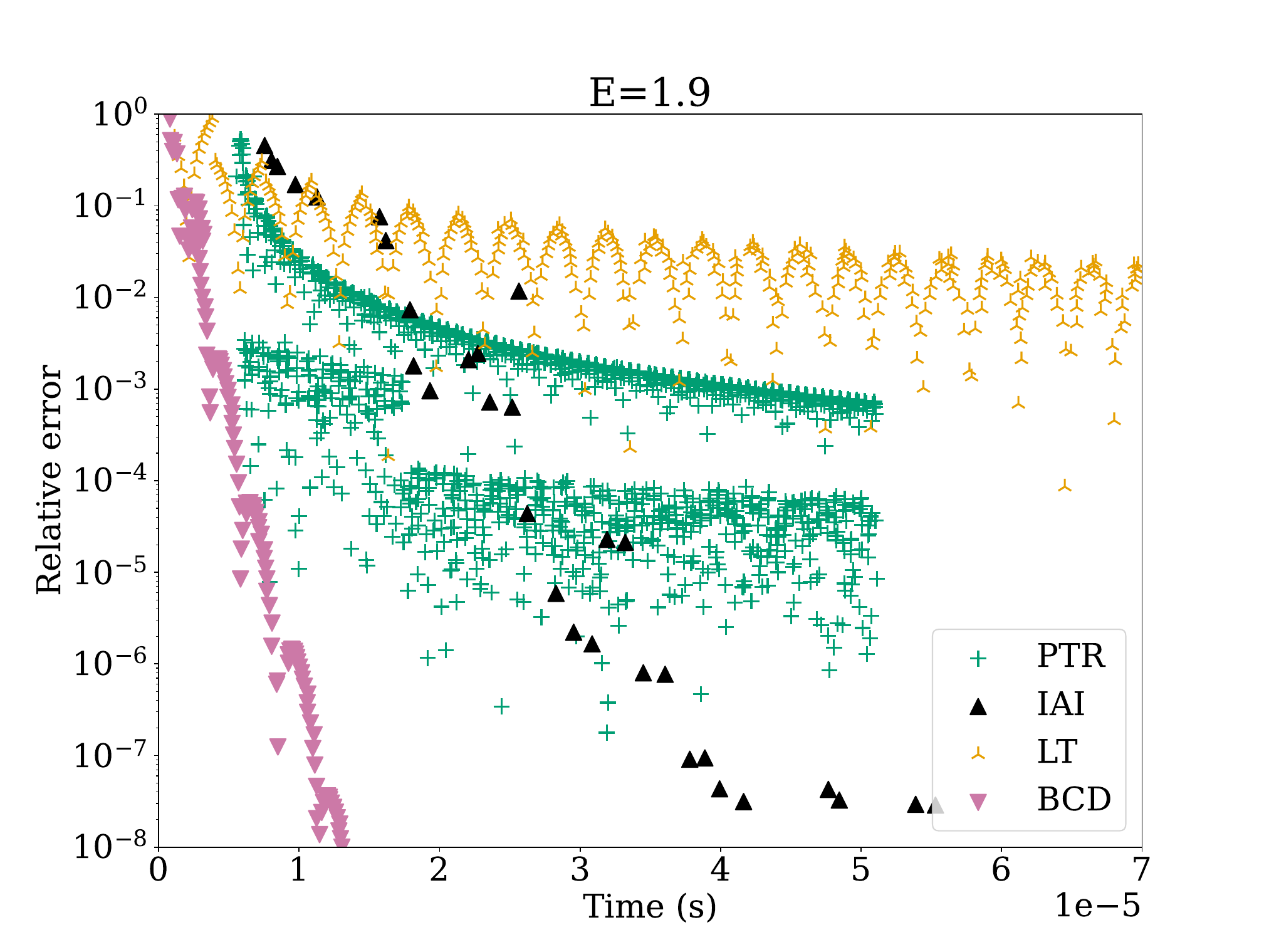}}

	\subfloat[Hard]{\includegraphics[scale=0.24]{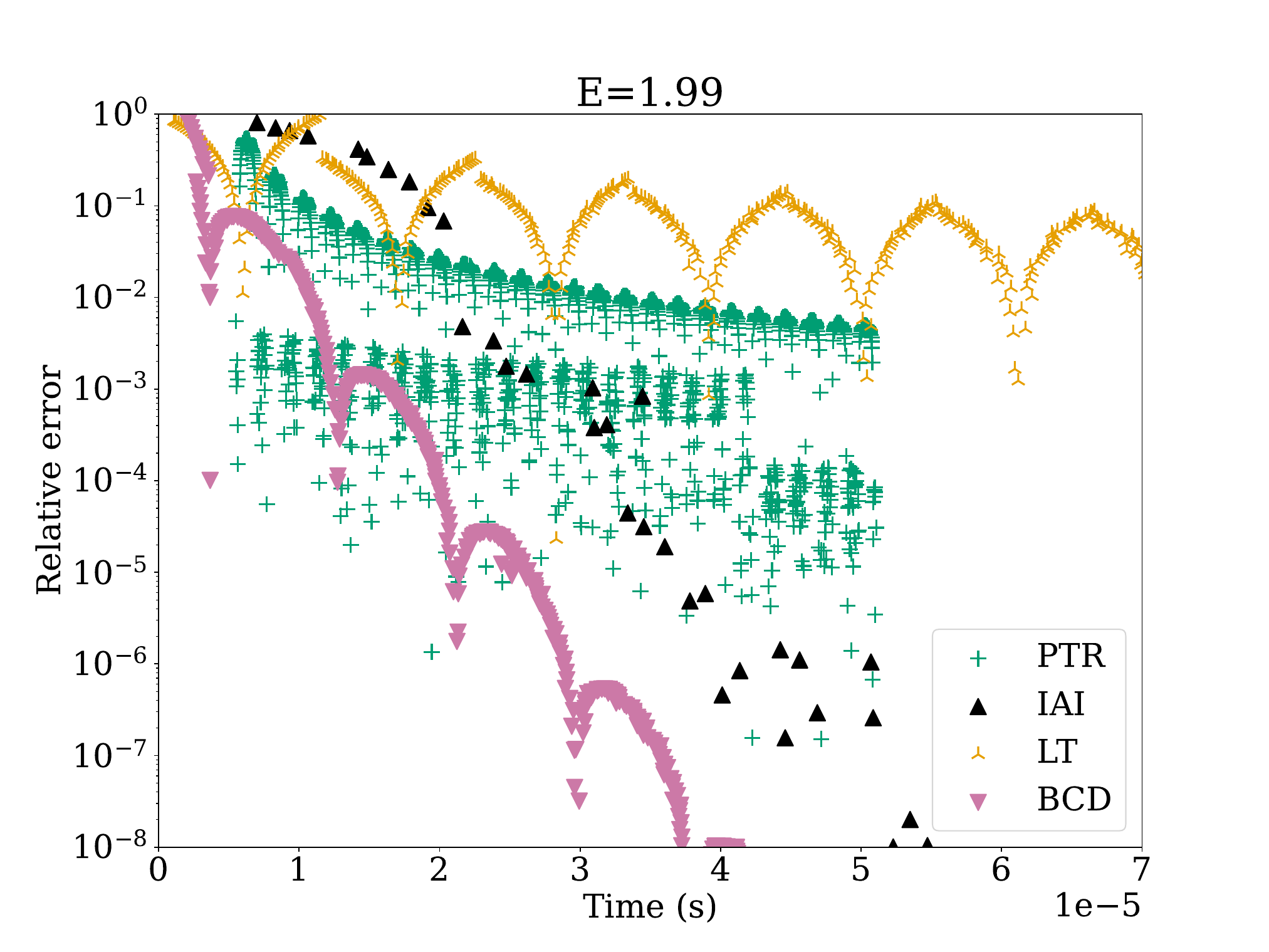}}
	\caption{Error as a function of computation time for the monatomic chain in 1D.}
	\label{fig:monatomictime}
\end{figure}

We summarize the results in~\autoref{fig:monatomictime}, where we respectively plotted the
relative error as a function of $N$ and of the computation time.

Away from van Hove singularities, the LT and  the PTR methods achieve errors no smaller than $10^{-6}$, whereas the IAI and the BCD attain machine precision. Near a van Hove, the LT’s accuracy drops to $10^{-3}$, while the IAI and the BCD retain full precision.
For all energies, the IAI converges to the smeared DOS, so the remaining error reflects the difference between smeared and non-smeared DOS. Notably, the BCD method attains machine precision near van Hove singularities, despite the DOS integrand losing the analytic structure the method exploits.

Regarding efficiency, far from van Hove singularities ($E=0.0$ eV), all methods yield comparable performance for relative errors above $10^{-3}$ whereas the BCD achieves superior accuracy for stricter tolerances ($<10^{-3}$).
For the medium case ($E=1.9$ eV), the PTR, the LT, and the BCD methods are most efficient for errors above $10^{-3}$, but BCD becomes the method of choice for higher accuracy.
In the vicinity of the van Hove singularity (hard case, $E=1.99$ eV), the optimal method remains the BCD for errors smaller than $10^{-1}$.

These results indicate that the PTR algorithm performs best for coarse errors across all energies, the LT is efficient for moderate errors away from van Hove singularities, the IAI excels at intermediate-to-high precision, and the BCD consistently provides the highest-precision computation.

\subsection{Two-band model of graphene}
We now consider the case of graphene, in the nearest-neighbour model
with hopping constant $t$ and lattice constant $a$ on an hexagonal
lattice with primitive lattice vectors $a_1 =\frac{a}{2} \left(3,\sqrt{3}\right)$
and $a_2 =\frac{a}{2} \left(3,-\sqrt{3}\right)$. This yields a $2\times 2$ Hamiltonian
with opposite eigenvalues.
Its reciprocal space Hamiltonian is

\begin{equation}
	H_{\mathbf{k}} = \begin{pmatrix}
		0                         & -t(1+e^{ik_x}+e^{ik_y}) \\
		-t(1+e^{-ik_x}+e^{-ik_y}) & 0
	\end{pmatrix},
\end{equation}
which eigenvalues are (for $\mathbf{k}=(k_x,k_y)$)

\begin{equation}
	\varepsilon_{\pm\mathbf{k}} =\pm \sqrt{3+2\cos\left( \sqrt{3}k_y a \right)+4 \cos\left( \frac{\sqrt{3}k_y a}{2} \right)\cos\left( \frac{3k_x a}{2} \right)}.
\end{equation}
with the Brillouin zone $\mathcal{B}$ an hexagonal cell centered around
$\Gamma=\mathbf{0}$ and spanned by the reciprocal-lattice vectors
$b_1 = \frac{2\pi}{3a}\left(1,\sqrt{3}\right)$, $b_2 = \frac{2\pi}{3a}\left(1,-\sqrt{3}\right)$.

The exact DOS is computed from the analytical band structure~\cite{castroneto2009electronic} using~\eqref{eq:dosgrad} and is given by
\begin{equation}
	D^{\mathrm{(C)_6}}_{\mathrm{ref}}(E)=\left\{\begin{matrix}
		\frac{2E}{t^2\pi^2}\frac{1}{\sqrt{F\left(E/t\right)}}K\left(\frac{4E/t}{F(E/t)}\right), & \ 0<|E|<t,          \\
		\frac{2E}{t^2\pi^2}\frac{1}{\sqrt{4E/t}}K\left(\frac{F(E/t)}{4E/t}\right),              & \ t<|E|<3t,         \\
		0,                                                                                      & \mathrm{elsewhere}.
	\end{matrix}\right.
\end{equation}
where $K$ is the complete elliptic integral of the first kind and $F(x)=(1+x)^2-(x^2-1)/4$.

\begin{figure}[ht!]
	\centering
	\includegraphics[scale=0.29]{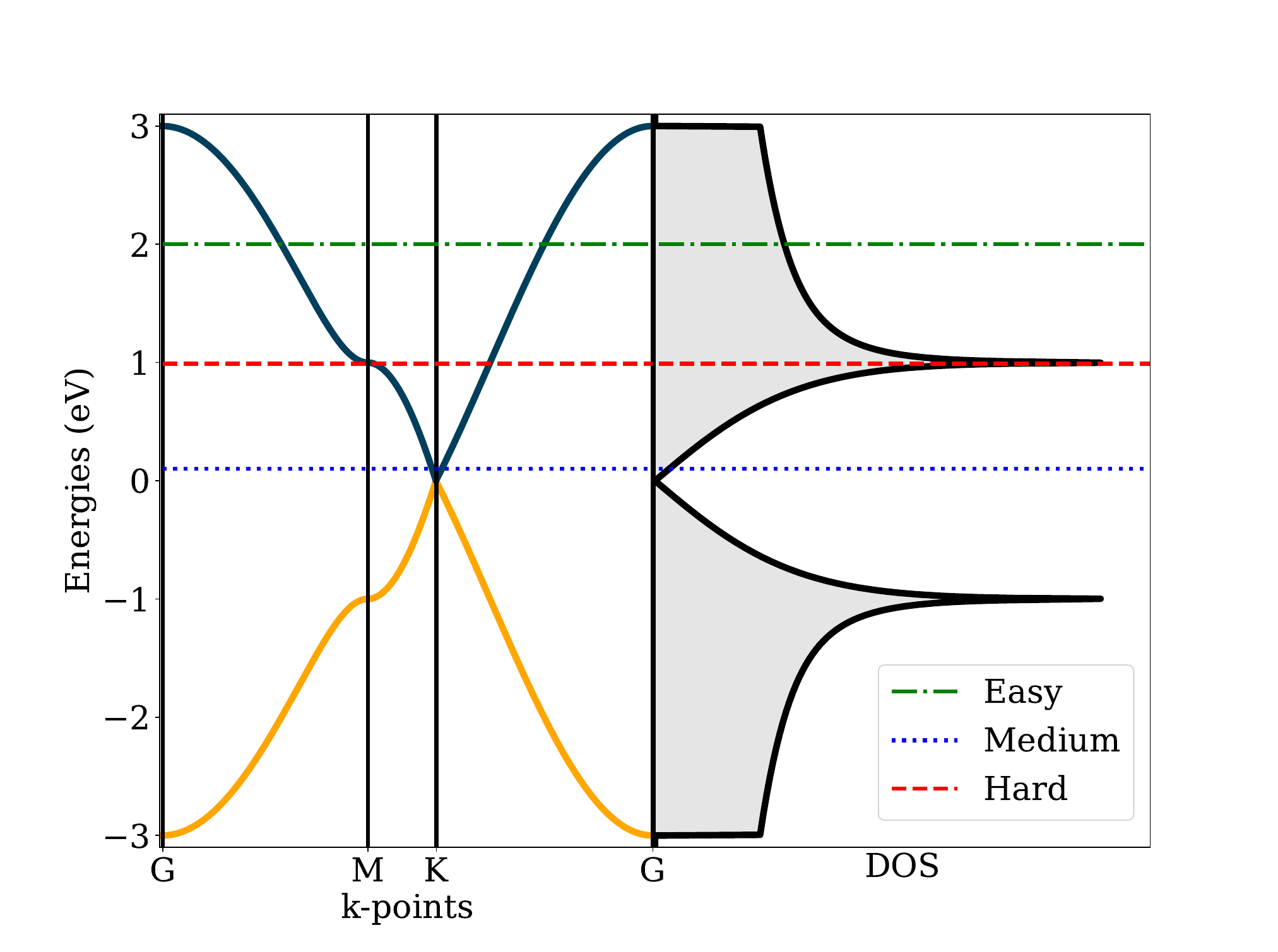}
	\caption{Bands and DOS of the graphene with hopping constant $t=1$.}\label{fig:dosgraphene}
\end{figure}

The DOS exhibits two saddle-point van Hove singularities at energies
$E=\pm t$ (with a logarithmic blow-up in the DOS), two extrema van Hove
singularities at points $E= \pm 3t$ (with a discontinuity),
and a Dirac point at energy $E=0$ (with a cusp). The bands and corresponding
DOS are shown in~\autoref{fig:dosgraphene}.

We benchmarked the numerical methods for three representative energies
depending on their relative distance to the van Hove
singularities: an easy case, in a smooth region ($E=2$), a medium case,
somewhat close to a Dirac point ($E=0.1$), and a hard case,
very close to a van Hove singularity ($E=0.99$).

\begin{figure}[ht!]
	\centering
	\subfloat[Easy]{\includegraphics[scale=0.24]{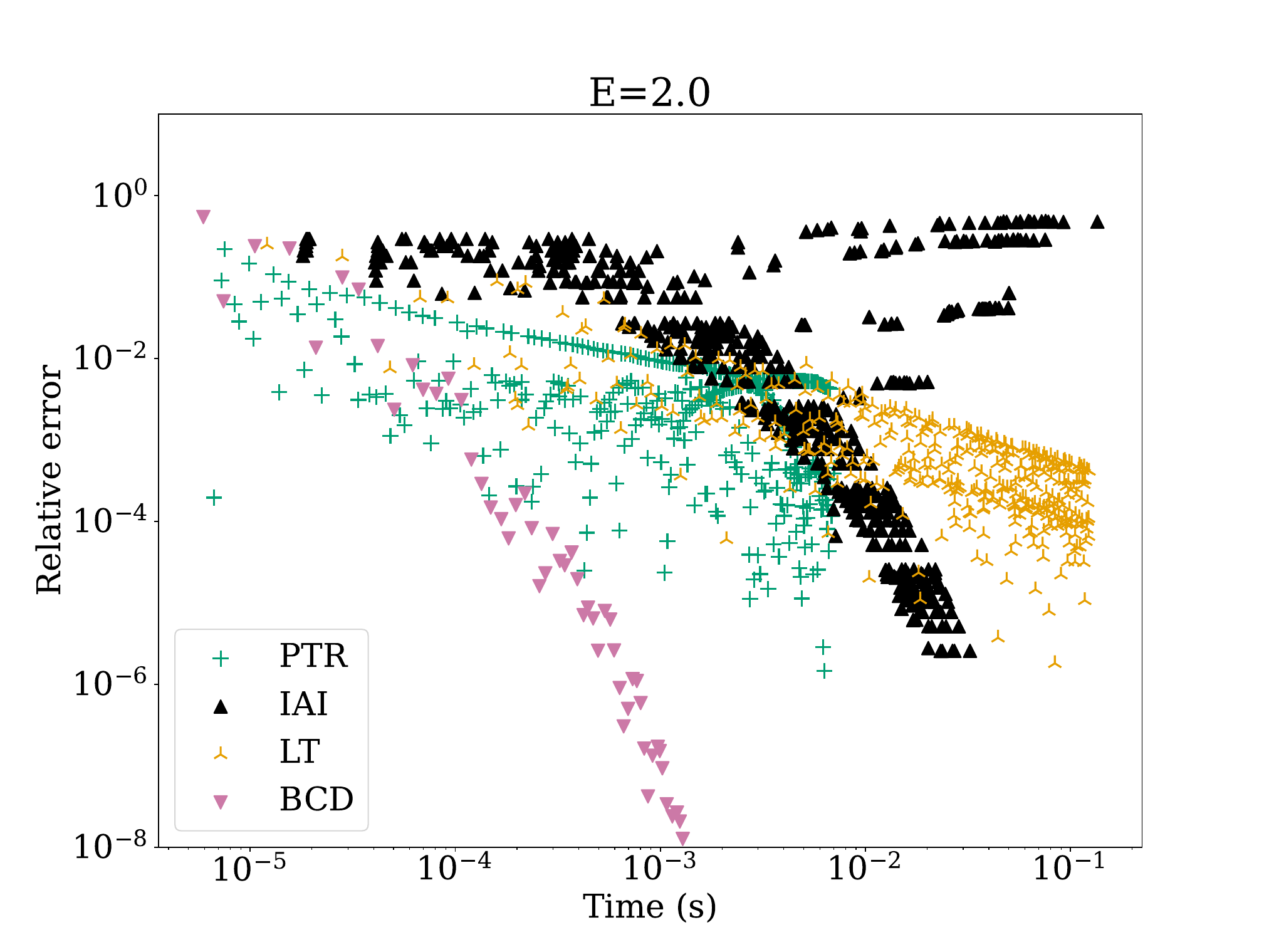}}
	\subfloat[Medium]{\includegraphics[scale=0.24]{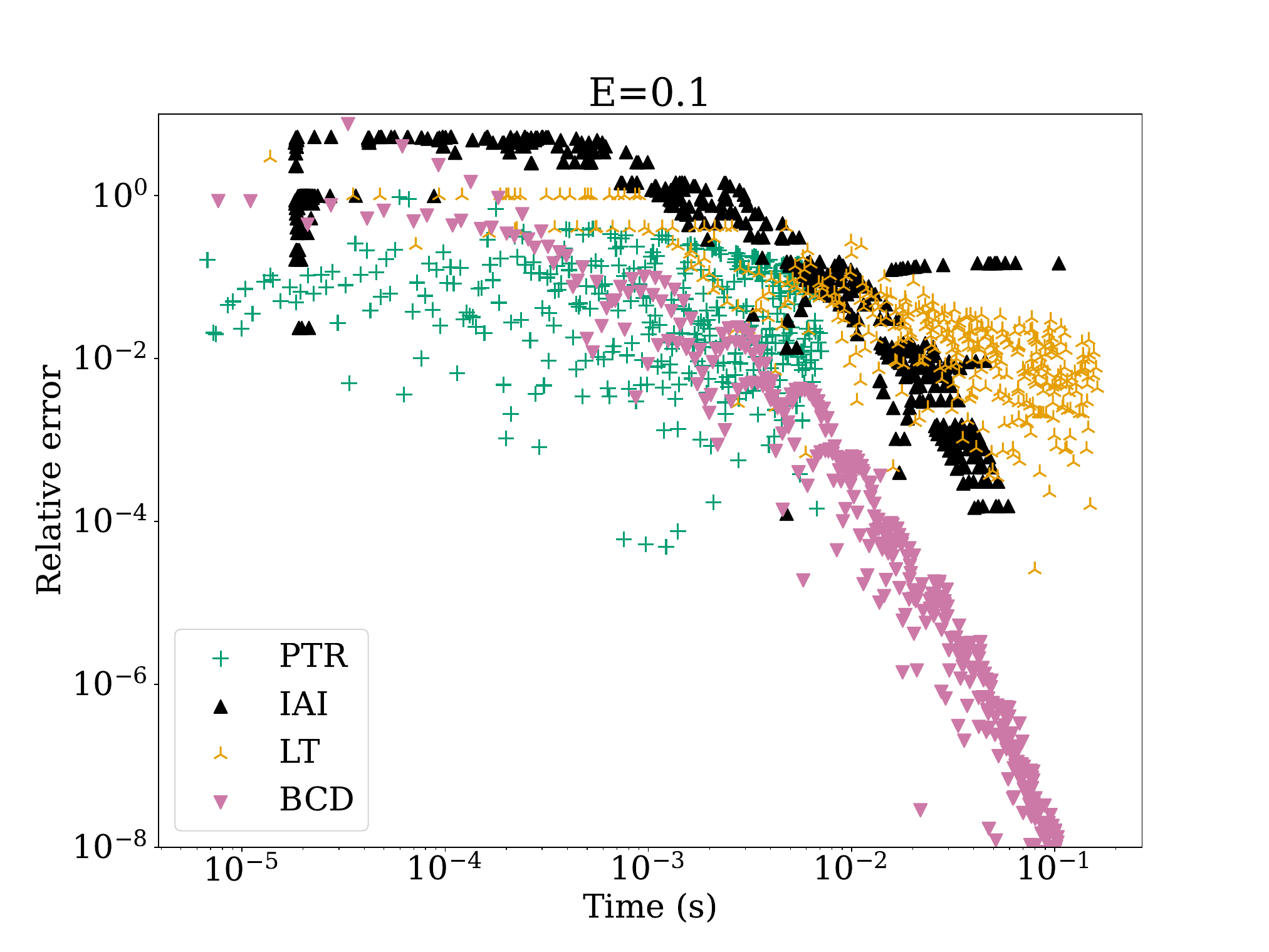}}

	\subfloat[Hard]{\includegraphics[scale=0.24]{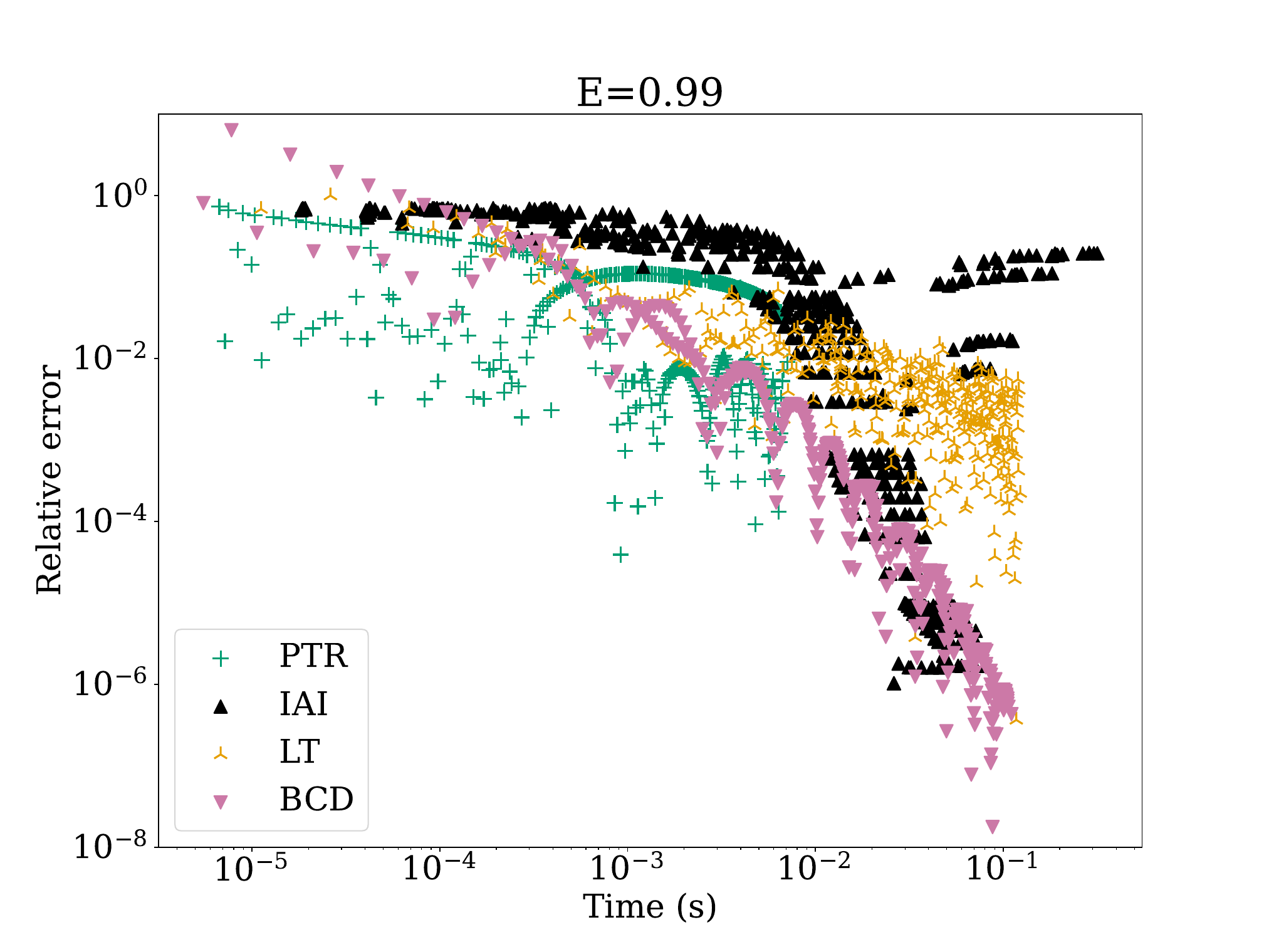}}
	\caption{Error as a function of computation time for graphene.}
	\label{fig:graphenetime}
\end{figure}

We summarize the results in~\autoref{fig:graphenetime}, where we plotted the
relative error as a function of the computation time.
Away from van Hove singularities (easy case, $E=2$ eV), all methods achieve errors below $10^{-5}$, with the BCD method reaching machine precision. Near the Dirac point (medium case, $E=0.1$ eV), the PTR, the LT, and the IAI plateau at errors around $10^{-3}$, while BCD attains errors below $10^{-8}$. In the vicinity of the van Hove singularity (hard case, $E=0.99$ eV), all methods exhibit increased errors: the PTR and the LT algorithms reach approximately $10^{-3}$, the IAI attains $10^{-6}$, and the BCD achieves $10^{-7}$.

Regarding efficiency, for the easy case ($E=2$ eV), the PTR is the most efficient method for errors above $10^{-3}$, while the BCD is preferable for higher accuracy. For the medium case ($E=0.1$ eV), the PTR remains optimal for errors above $10^{-4}$, with the BCD outperforming for smaller errors. Finally, near the van Hove singularity (hard case, $E=0.99$ eV), the PTR is most efficient for errors above $10^{-2}$, the IAI and the BCD methods are optimal for errors between $10^{-3}$ and $10^{-6}$, and the BCD is superior for errors below $10^{-6}$.

These results indicate that the PTR is best-fitted for low-accuracy, the IAI excels for intermediate accuracy near singularities, and the BCD consistently provides the highest precision.

\subsection{Strontium Vanadate (\texorpdfstring{$\mathrm{SrVO_3}$}{SrVO3})}

We study the DOS of the strontium vanadate ($\mathrm{SrVO_3}$), following the
example of~\cite{kaye2023automatic}. We focus on the three $\mathrm{t_2g-3d}$
orbitals on a primitive cubic lattice, which yields a $3\times 3$
Hamiltonian with degenerate eigenvalues. This Hamiltonian is computed
using Wannierization on a $11\times 11 \times 11$ k-point mesh of the Brillouin
zone $\mathcal{B}= 2\pi\mathbb{R}^3/\mathbb{Z}^3$. Its associated chemical potential is
$\varepsilon_F = 12.39 \ \mathrm{eV}$.

\begin{figure}[ht!]
	\centering
	\includegraphics[scale=0.29]{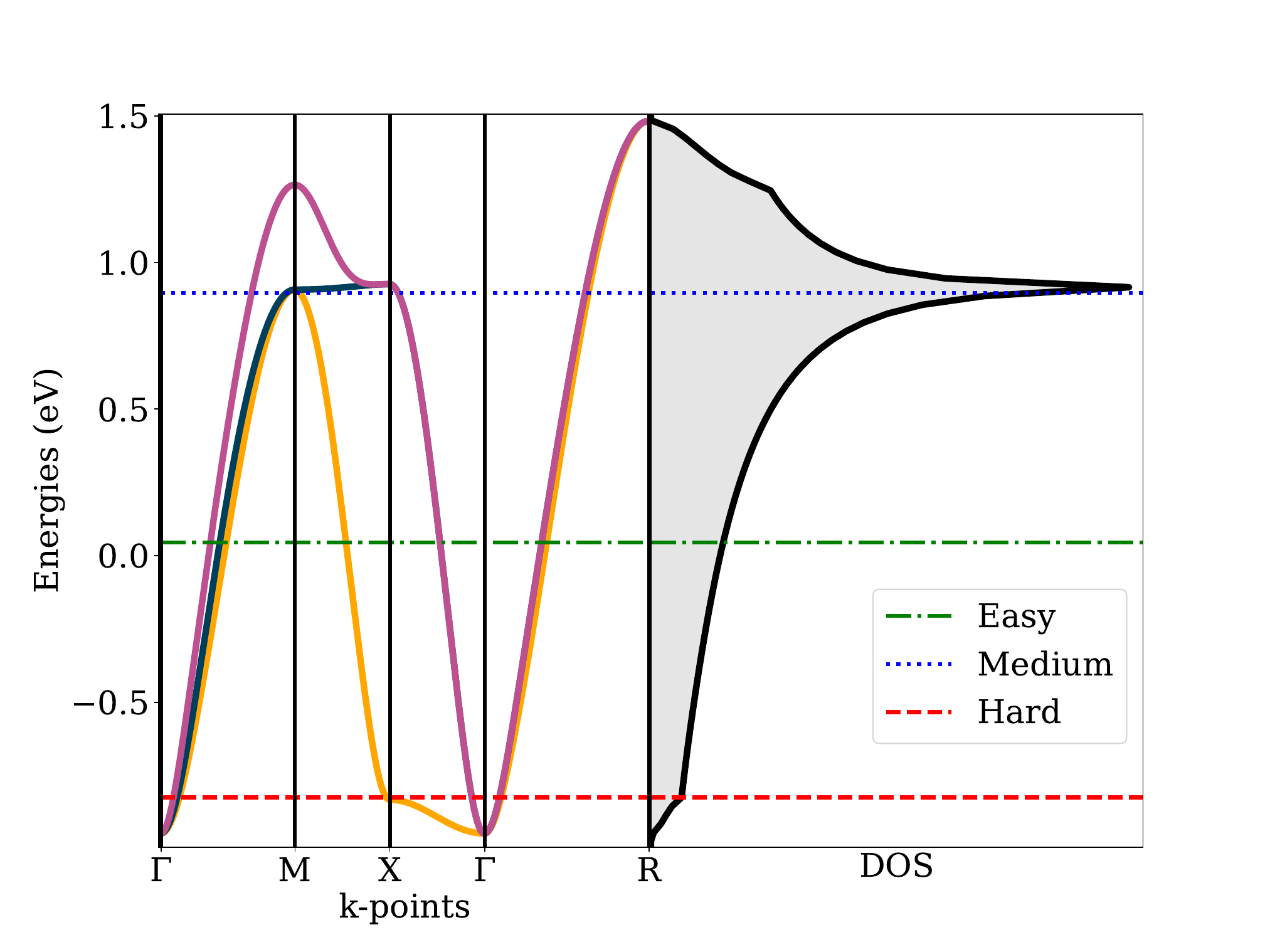}
	\caption{Bands and DOS of $\mathrm{SrVO_3}$.}\label{fig:dossrvo3}
\end{figure}

The DOS exhibits three van Hove singularities at energies
$E=-0.82 \text{ eV}$, $E=0.92 \text{ eV}$ and
$E=1.23 \text{ eV}$. As is typical for three-dimensional
system~\cite{vanhove1953occurrence}, the DOS does not diverge at van Hove
singularities and only exhibits cusps. The reference DOS against which all relative errors are measured was generated with the BCD algorithm on the irreducible, symmetry-reduced Brillouin zone, with $N=501$. Convergence was confirmed by an independent calculation using the IAI method with a small smearing parameter $\eta$. The bands and corresponding DOS are shown in~\autoref{fig:dossrvo3}.
We benchmarked the numerical methods for three representative energies
depending on their relative distance to the van Hove
singularities: an easy case, in a smooth region ($E=0.05 \text{ eV}$), a medium case, close to a van Hove singularity ($E=0.90 \text{ eV}$), and a hard case,
somewhat close to a small cusp ($E=-0.84 \text{ eV}$).

\begin{figure}[ht!]
	\centering
	\subfloat[Easy]{\includegraphics[scale=0.24]{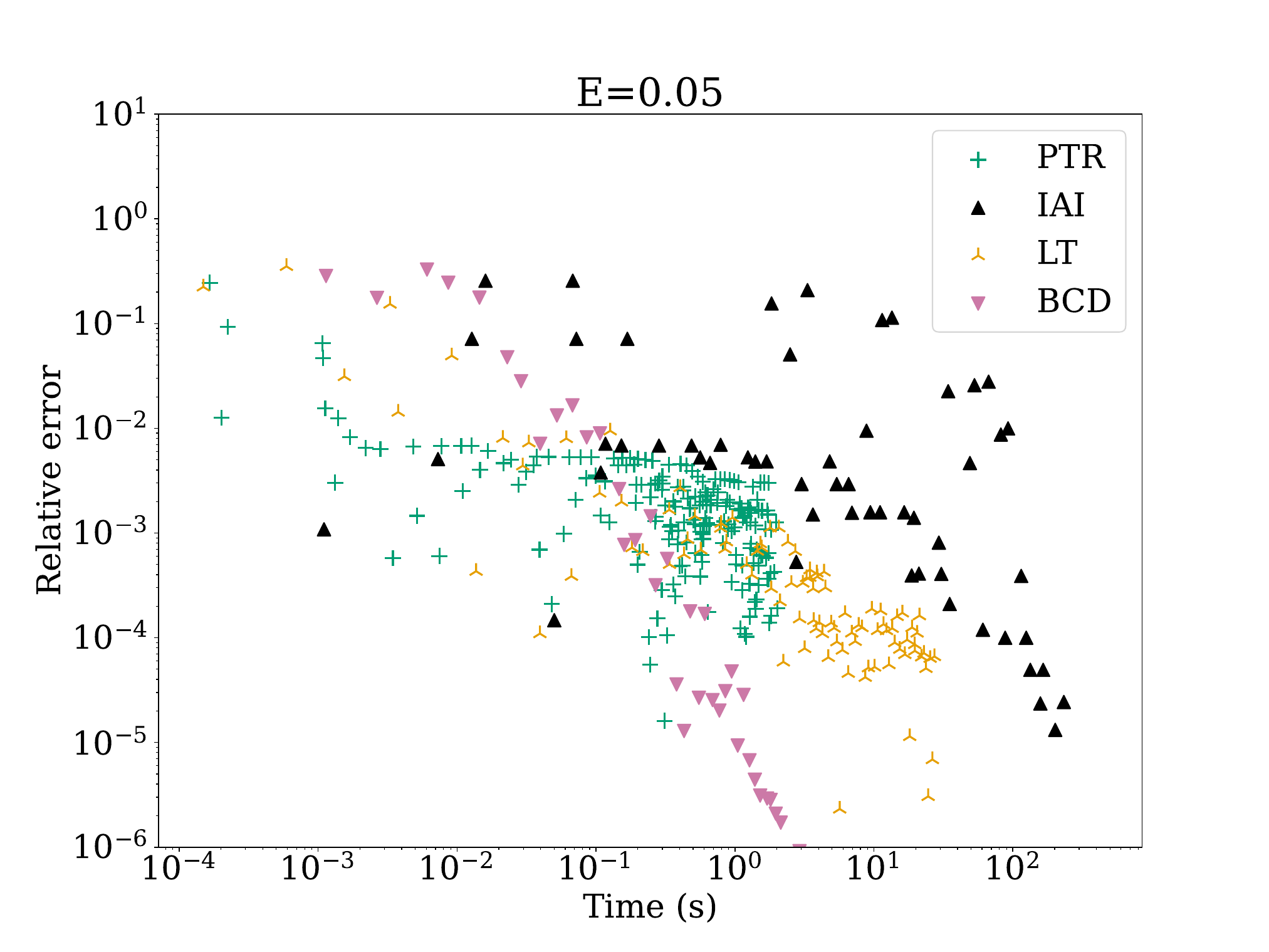}}
	\subfloat[Medium]{\includegraphics[scale=0.24]{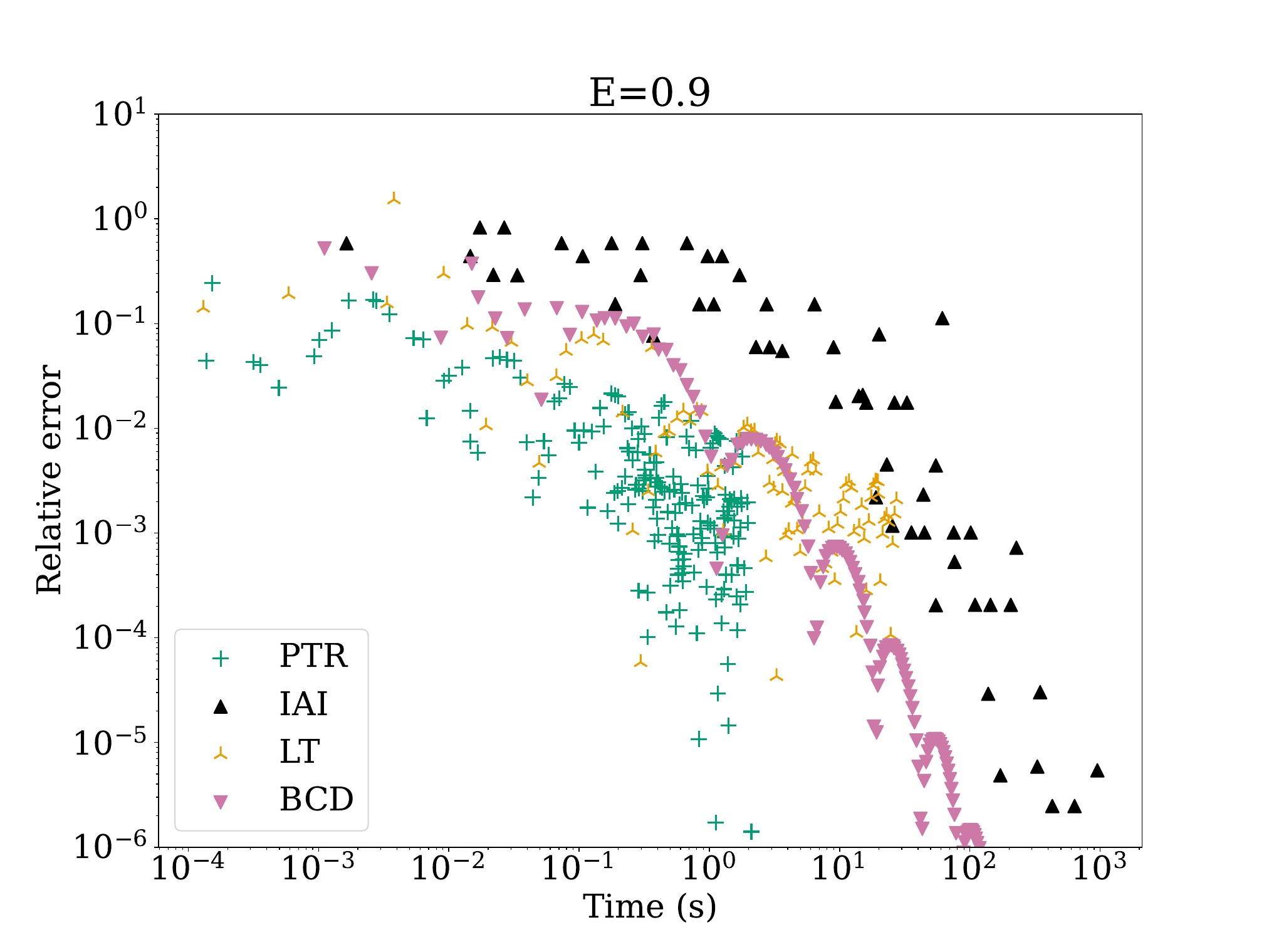}}

	\subfloat[Hard]{\includegraphics[scale=0.24]{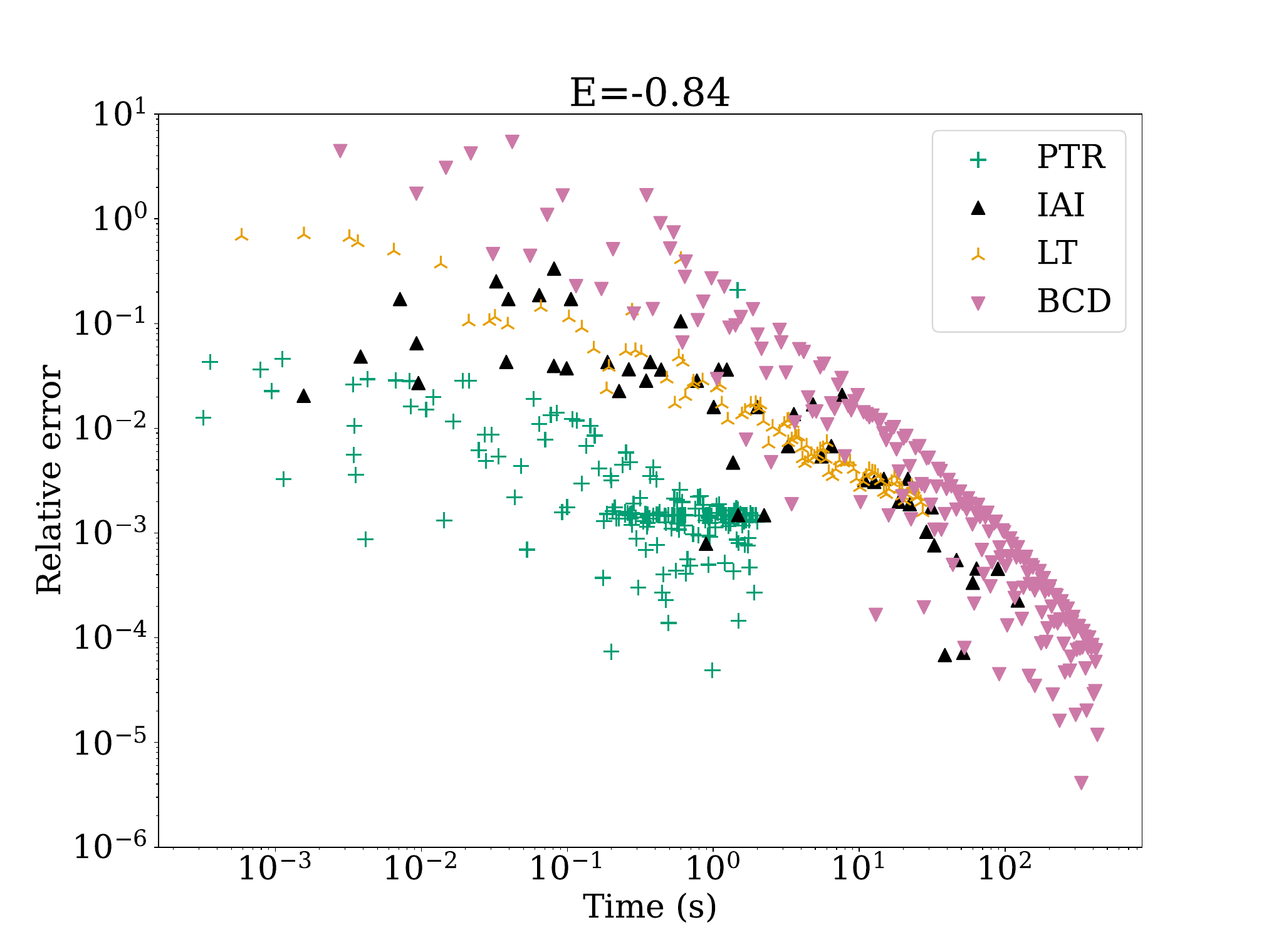}}
	\caption{Error as a function of computation time for $\mathrm{SrVO_3}$.}
	\label{fig:srvo3time}
\end{figure}

We summarize the results in~\autoref{fig:srvo3time}, where we plotted the
relative error as a function of the computation time.
Away from van Hove singularities (easy case, $E=0.05$ eV), all methods achieve errors below $10^{-5}$, with the BCD reaching machine precision. Near a van Hove singularity (medium case, $E=0.90$ eV), the PTR and LT algorithms plateau at errors around $10^{-3}$, the IAI attains $10^{-5}$, while the BCD method achieves errors below $10^{-7}$. In the vicinity of a small cusp (hard case, $E=-0.84$ eV), the PTR and the LT reach approximately $10^{-3}$, the IAI and BCD algorithms achieve $10^{-5}$ relative error.

Regarding efficiency, for the easy case ($E=0.05$ eV), all methods are equivalently performing for relative error above $10^{-5}$, whereas the BCD outperforms the other methods for higher precision. For the medium case ($E=0.90$ eV), the PTR is the best method for errors above $10^{-4}$, but for smaller errors BCD is the best-fitted method. Finally, near the cusp (hard case, $E=-0.84$ eV), the PTR outperforms the other methods for error above $5\times 10^{-4}$, whereas the BCD and the IAI methods yield similar results for smaller relative error down to $10^{-5}$.

\section{Smeared DOS Results}\label{section:smeared_dos}

\begin{figure}[ht!]
	\centering
	\subfloat[Easy]{\includegraphics[scale=0.24]{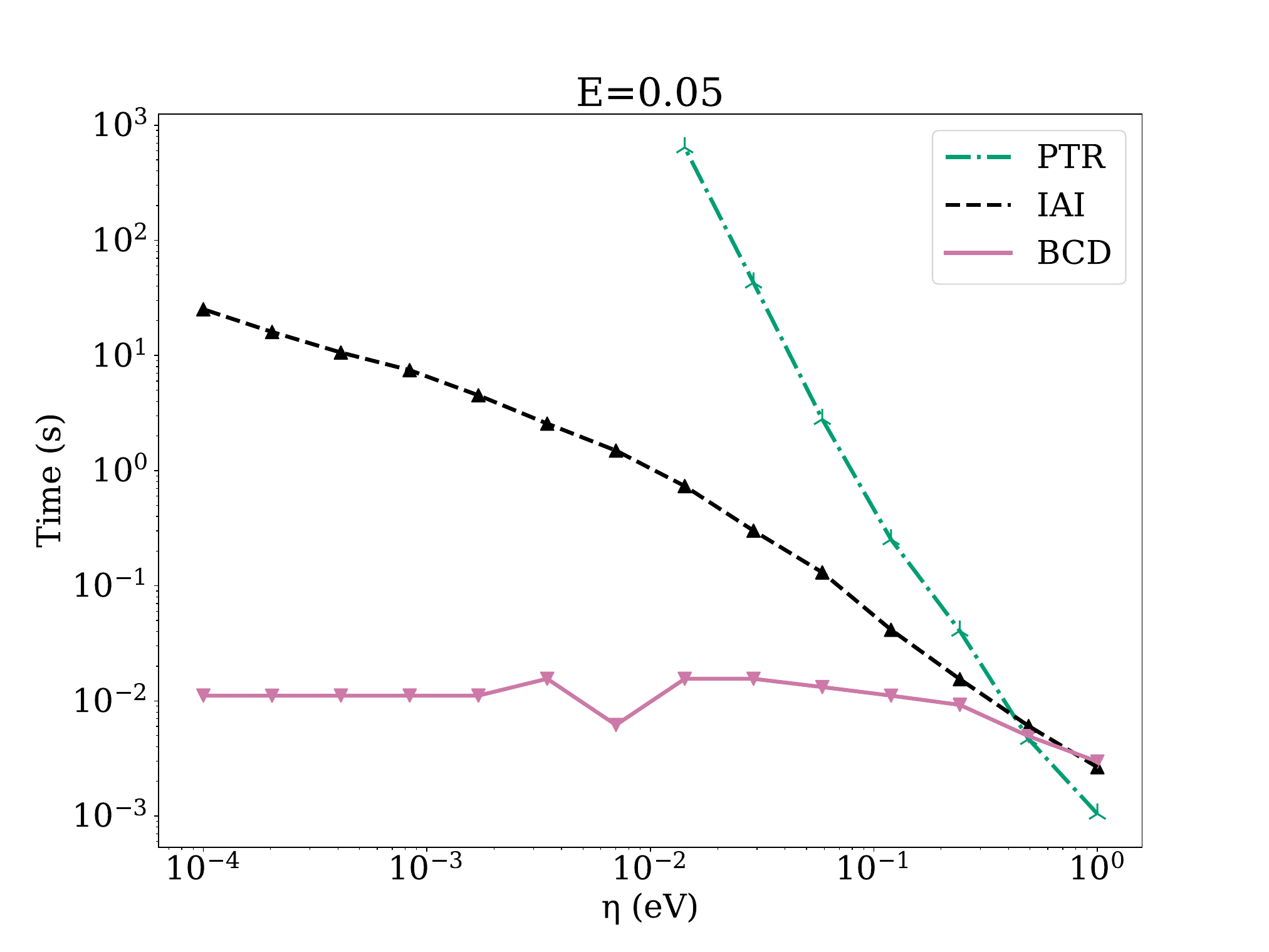}}
	\subfloat[Medium]{\includegraphics[scale=0.24]{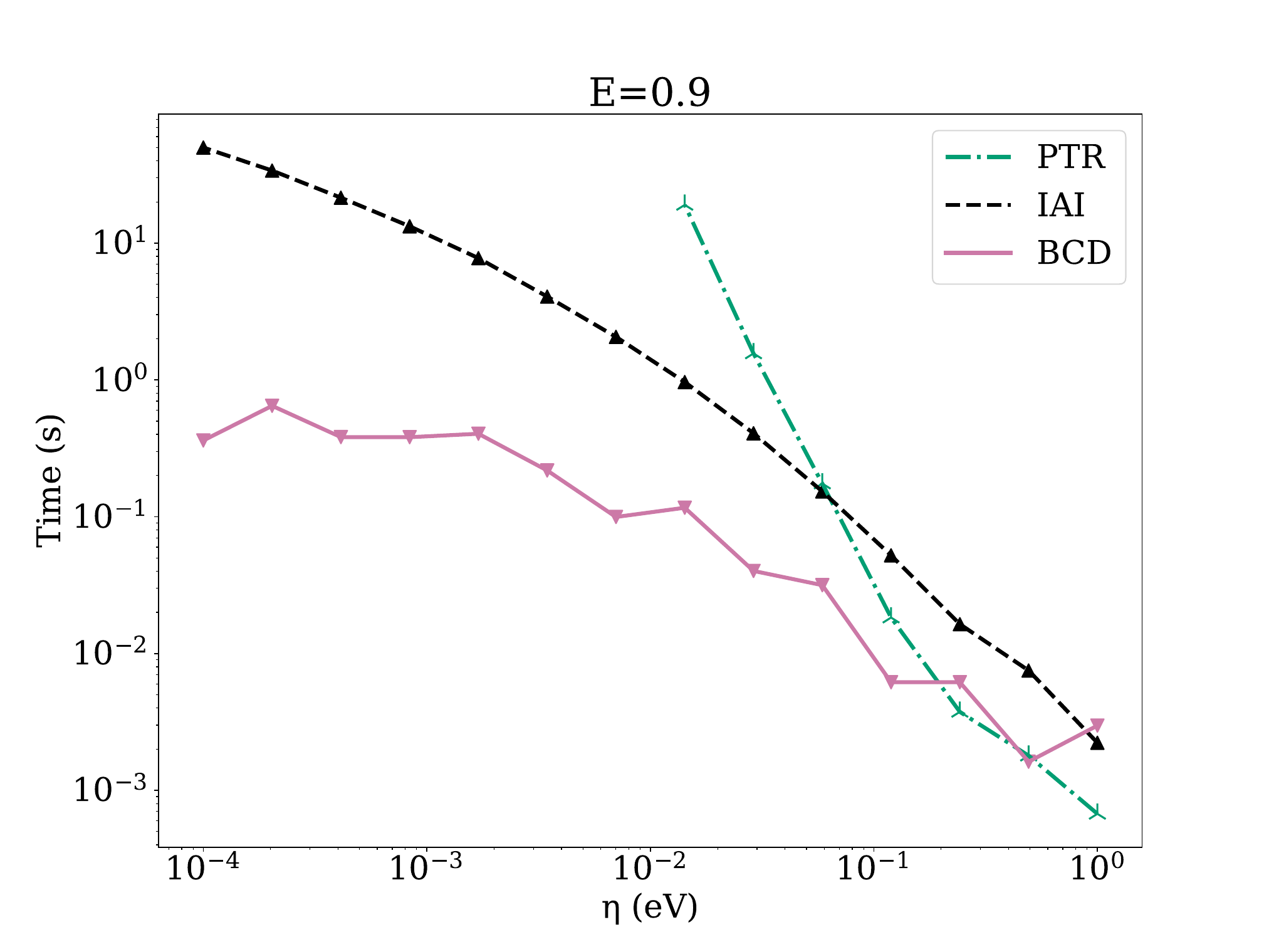}}

	\subfloat[Hard]{\includegraphics[scale=0.24]{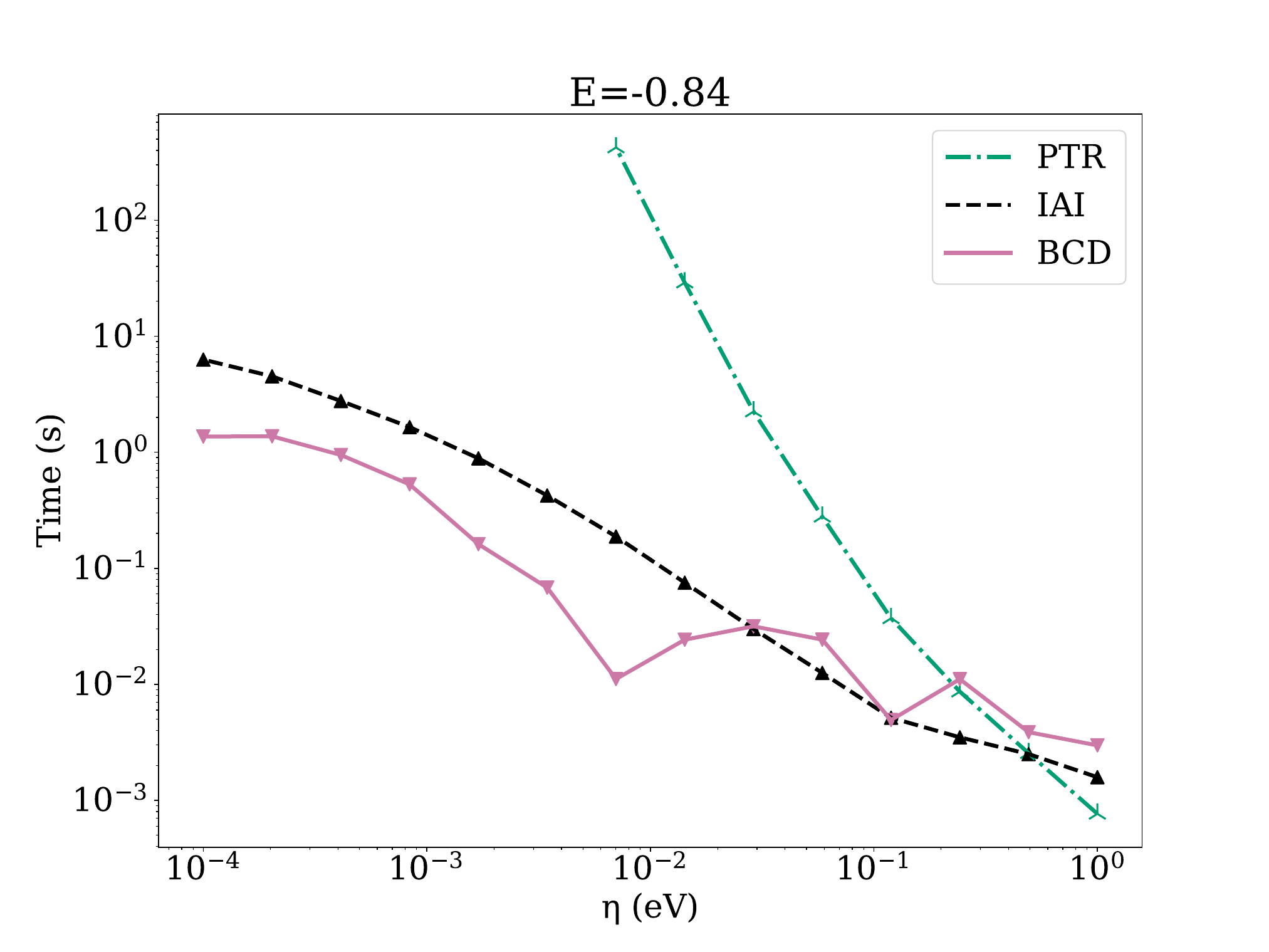}}
	\caption{Time required to reach the absolute error $\epsilon = 10^{-5}$ as a function of $\eta$.}
	\label{fig:errorsmearing}
\end{figure}

Here we report the computational cost required to reach a fixed accuracy $\epsilon =10^{-5}$ for the smeared density of states (DOS) of $\mathrm{SrVO_3}$. The dependence on the broadening parameter $\eta$ is shown in~\autoref{fig:errorsmearing} for three representative energies corresponding to increasing numerical difficulty.

As expected, the asymptotic scaling behaviors are recovered. The PTR exhibits a cost scaling as $\mathcal{O}(\eta^{-1})$, while the IAI follows $\mathcal{O}(\log(\eta^{-1}))$.
In contrast, the BCD shows an approximately constant computational cost, $\mathcal{O}(1)$ when the energy lies sufficiently far from van Hove singularities. However, a clear $\eta$-dependence emerges in their vicinity due to integrand non-analyticity near singularities.

From a practical perspective, ~\autoref{fig:errorsmearing} shows that PTR is the most efficient method for large broadening, typically $\eta \gtrsim 10^{-1}$.
For smaller $\eta$, the BCD becomes the most efficient approach when the energy is far from singularities, whereas the IAI appears to be a good alternative in their vicinity.
\section{Conclusion}\label{section:conclusion}
A central contribution of this work is a detailed analysis of the
domain of validity of the BCD method. We have identified and
characterised its systematic failures at non-van Hove energies
associated with symmetry-protected band crossings. While it is
unusable for artificial high-symmetry systems possessing a large
number of band crossings such as the free electron gas, it works
without issue for the realistic materials we have studied. To detect
these failures in practice, we have proposed monitoring the sign of
the imaginary part of the eigenvalues of the deformed Hamiltonian
$H_{\mathbf{k}+i\mathbf{h}(\mathbf{k})}$ as a diagnostic.

Our results establish clear performance regimes for each method. For
non-smeared DOS calculations, the PTR algorithm is optimal for low-precision
requirements (relative error above $10^{-2}$), while the  BCD is more
efficient in the high-accuracy regime and not too close to van Hove
singularities. The IAI offers a reliable compromise for intermediate
precision, particularly effective near singularities. For smeared DOS,
the PTR remains optimal for large smearing parameters ($\eta>0.1$), while
the IAI and the BCD become superior for smaller $\eta$, with the BCD providing the highest efficiency away from problematic crossings. The LT method
does not appear to be competitive except for low-accuracy estimations.

Apart from the above-mentioned deficiency of the BCD in highly-symmetric
systems, let us note that the BCD, like the IAI, works natively
energy-by-energy and computations done for a given energy can only be
reused in a limited fashion to other nearby energies. Future work
could explore extensions of the BCD to handle symmetric crossings more
robustly as well as work for several energies at once, for instance by
using component-dependent deformations. The handling of van Hove
energies (and its neighborhoods) also remains an issue.

\section*{Acknowledgments}
We are grateful to Pierre Amenoagbadji, Jason Kaye and Lorenzo Van
Mu\~noz for stimulating discussions.
\sloppy
\printbibliography
\fussy
\appendix
\section{Numerical computation of derivatives}
\label{divdif}
We explain here how to compute the derivative of
\begin{align*}
	\mathbf{h}(\mathbf k) & =-\alpha \left(\mathrm{Tr}\left(\nabla_{\mathbf k} H_\mathbf{k}\chi\left(\frac{H_\mathbf{k}-E}{\Delta E}\right)\right)\right).
\end{align*}
For each component $j =1, \dots, d$ of $\mathbf{h}$, we need to compute
the derivatives of
\begin{equation}
	f(\mathbf{k}) = \mathrm{Tr}\left(\partial_{j}H_{\mathbf{k}}\chi\left(H_{\mathbf{k}}\right)\right),
\end{equation}
where we have set $\Delta E, -\alpha$ to $1$ and $E$ to $0$ for notational simplicity. We have, for all
$i=1, \dots, d$,
\begin{equation}
	\partial_i f(\mathbf{k}) = \mathrm{Tr}(\partial_i\partial_j H_{\mathbf{k}} \chi(H_{\textbf{k}})) + \mathrm{Tr}\left( \partial_j H_{\textbf{k}} d\chi\left(H_{\textbf{k}} \right)\cdot \partial_i H_{\textbf{k}}\right).
\end{equation}
This function remains smooth under the assumptions given in
\eqref{eq:cond1} and~\eqref{eq:cond2}. The first term is
straightforward to compute. Using perturbation theory, the second term
can be expressed in a basis of orthonormal eigenvectors as~\cite{higham2008functions}
\begin{equation}
	\mathrm{Tr}\left( \partial_j H_{\textbf{k}} d\chi\left(H_{\textbf{k}} \right)\cdot \partial_i H_{\textbf{k}}\right) = \sum_{m,n} \frac{\chi(\varepsilon_n)-\chi(\varepsilon_m)}{\varepsilon_n-\varepsilon_m} \langle   \psi_n|\partial_j H_{\textbf{k}}|\psi_m\rangle\langle \psi_m|\partial_i H_{\textbf{k}}|\psi_n\rangle,
\end{equation}
where we use the convention
\begin{equation}
	\frac{\chi(\varepsilon_n)-\chi(\varepsilon_n)}{\varepsilon_n-\varepsilon_n} = \chi'(\varepsilon_n).
\end{equation}

Numerically, the naive evaluation of the divided difference
$(\chi(\varepsilon_n)-\chi(\varepsilon_m))/(\varepsilon_n-\varepsilon_m)$
can introduce numerical errors when several bands are close together,
as small inaccuracies (of the order of the machine epsilon) in the
computation of $\chi(\varepsilon_n)-\chi(\varepsilon_m)$ will yield
potentially large inaccuracies in the computation of
$\tfrac{\chi(\varepsilon_n)-\chi(\varepsilon_n)}{\varepsilon_n-\varepsilon_n}$.
There are general ways around this problem \cite{quan2025matrixfuns},
but for particular cutoff functions one can reformulate the expression
to give numerically accurate results. For instance, for the gaussian,
\begin{equation}
	\frac{\exp(-x^2)-\exp(-y^2)}{x-y} =\exp(-y^2)\frac{\mathrm{\tt expm1}(-(x-y)(x+y))}{x-y},
\end{equation}
where \texttt{\tt expm1} is a widely available function that computes
$x \mapsto e^{x} -1$ in a numerically robust way for small positive values. More
generally, such algebraic tricks can always be performed
\cite{reps2003computational}.
\end{document}